\documentclass[10.5pt, reqno,twoside,onecolumn,final]{siamltex}
\usepackage{amscd}
\usepackage{amssymb}
\usepackage{amsmath,xcolor}
\usepackage{mathrsfs}
\usepackage{graphicx}
\usepackage{calc}%
\usepackage{multirow}
\usepackage{placeins}
\usepackage{cite}
\usepackage{mathptmx}  

\newcommand{\m}{\mathbf m }

\newcommand{\vv}{\mathbf v}
\newcommand{\x}{\mathbf x}
\newcommand{\y}{\mathbf y}
\newcommand{\kk}{\mathbf k}
\newcommand{\bfa}{\mathbf a}
\newcommand{\alp}{\boldsymbol \alpha}

\newcommand{\diver}{{\rm {div} \,\, }}

\newcommand{\R} {\mathbb R}

\newcommand{\nablad} {\nabla\cdot}

\newcommand{\intd}[1]{\left( #1 \right)}

\newcommand{\ddt}{\frac d{dt} }

\newcommand{\Dt}{\Delta t}

\def\eqdef{\stackrel{\rm def}{=}}

\def\beq{\begin{equation}}
\def\eeq{\end{equation}}
\def\beqs{\begin{equation*}}
\def\eeqs{\end{equation*}}

\newcommand{\esssup}{\mathop{\mathrm{ess\,sup}}}

\newcommand{\varep}{\varepsilon}
\newcommand{\brho}{\bar\rho}
\newcommand{\norm}[1]{\left\| #1 \right\|  }

\newcommand{\indic}{{\mathbf 1}}

\numberwithin{equation}{section}

\def\myclearpage{}
\usepackage[left=1 in,top=1.25 in,right=1 in,bottom=1.25 in]{geometry}

\title{A Mixed Finite Element Approximation for Fluid Flows of Mixed Regimes in Porous Media}
\author{John Cummings, Matthew Hamilton  and  Thinh Kieu\footnotemark[1]  }
\date{today}

\begin{document}
\maketitle
\renewcommand{\thefootnote}{\fnsymbol{footnote}}

\footnotetext[1]{Department of Mathematics, University of North Georgia, Gainesville Campus, 3820 Mundy Mill Rd., Oakwood, GA 30566, U.S.A. ({\tt email: John.Cummings@ung.edu, MAHAMI7097@ung.edu, thinh.kieu@ung.edu}).}               
               
\begin{abstract} 
In this paper, we consider the complex flows when all three regimes pre-Darcy, Darcy and post-Darcy may be present in different portions of a same domain. 
We unify all three flow regimes under mathematics formulation.
We describe the flow of a single-phase fluid in $\R^d, d\ge 2$ by a nonlinear degenerate system of density and momentum. A mixed finite element method is proposed for the approximation of the solution of the above system.   
The stability of the approximations are proved; the error estimates are derived for the numerical approximations for both continuous and discrete time procedures.  
 The continuous dependence of numerical solutions on physical parameters are demonstrated. 
Experimental studies  are presented regarding convergence rates and showing the dependence of the solution on the physical parameters.

\end{abstract}
            
            



\pagestyle{myheadings}
\thispagestyle{plain}
\markboth{John Cummings, M. Hamilton and T. Kieu}{Mixed FEM for fluid flows of mixed regimes}
           
\myclearpage    
\section{Introduction}\label{intro}
Fluid flow in porous media plays an important role in a wide range of science and engineering applications, e.g water resources, geothermal system, chemical processes, gas and water purification, gas storage, oil extraction and petroleum engineering. Fluid flow in porous media is very complicated and are modeled by many different equations of various types due to the variety and complexity of the filtration matrix of the media in the path of the fluid. Based on the experimentally observed non‐Darcy behavior of seepage flows by Basak \cite{BASAK77}, Soni et al. \cite{SONI1978231}, Kundu et al. \cite{KUNDU2016278} classified flow regimes into three main flow regimes: pre‐Darcy flow (i.e. pre-linear, non-Darcy), Darcy flow (linear), and post‐Darcy flow (i.e. post-linear, non-Darcy). There is a general consensus that the Darcy regime is valid as long as the Reynolds number ($Re$) in the range of characteristic value between $1$ and $10$, see in \cite{BearBook}. When the Reynolds number is high ($Re>10$), there is a deviation from the Darcy law and Forchheimer's equations are usually used to account for it \cite{ForchheimerBook}, see also in \cite{Muskatbook,BearBook,NieldBook}. 
On the other end of the Reynolds number's range, when it is very small ($Re\to 0$),  the pre-Darcy regime is observed but not fully understood, although it contributes to unexpected oil extraction, and improved oil recovery in petroleum reservoirs, see in \cite{Dudgeon85,SSHI2016,SoniIslamBasak78,Bloshanskaya2017} and references therein.

Most studies in porous media focus on the Darcy regime which is presented by the (linear) Darcy equation which is a linear relationship between the pressure gradient and the fluid velocity, see in \cite{VazquezPorousBook}.
Recently, the post-Darcy regime has been attracted attention with the (nonlinear)  mathematical and numerical modeling of Forchheimer flows, see in \cite{StraughanBook,ABHI1,HI2,HIK1,HIK2,HKP1,HK2,CHK1,CHK2, K2, K3} and references therein.
In contrast, the (nonlinear) pre-Darcy regime has not received much attention among the researchers and engineers.
Moreover, the three regimes are always treated separately.
This is due to the different natures of the models and the ranges of their applicability.
There are evidences in the literature that the flow of a fluid in porous media may present all three regimes in different unidentified portions of the confinement. 
The petroleum reservoir is an example of environment exhibiting all the three regimes: fast flow near wells and fractures is described by   Forchheimer equation, very slow flow far away from wells is described by the Pre-Darcy equation and moderate flow in the between is described by Darcy equation.       
Therefore, there is a need to unify the three regimes into one formulation and study the fluid as a whole. This paper aims at deriving admissible models for this unification and analyze their properties mathematically. 

We now start the investigation of different types of fluid flows in porous media.
Consider fluid flows with velocity $\vv\in \mathbb R^d$, $d\ge 2$ pressure $p\in \mathbb R$, and density $\rho\in[0,\infty)$.
Depending on the range of the Reynolds number, there are different groups of  equations to describe their dynamics.

Darcy's  law is the most popular equation to describe flow in porous media and fractures at moderate flow rates, when the
flow rate and the pressure gradient have a linear relationship. 
\beq\label{D}
\vv=- k(\x,t)\nabla p(\x,t), \text{ where $k(\x,t)>0$.} 
\eeq

However, as the flow rate (Reynolds’s number) decreases, there are Izbash-type  equations  that describe the pre-Darcy regime:
\beq\label{PD}
|\vv(\x,t)|^{-\alpha}\vv(\x,t)=- k(\x,t) \nabla p(\x,t)\, \text{for some constant power } \alpha\in (0,1)\text{ and  coefficient } k(\x,t)>0.
\eeq
For experimental values of  $\alpha$, see e.g.  \cite{SSHI2016,SoniIslamBasak78}.

When the flow rate (Reynolds’s number) is high, the following Forchheimer equations are widely used in studying post-Darcy flows. 

Forchheimer's two-term law
\beq\label{F2}
a(\x,t)\vv(\x,t)+b(\x,t)|\vv(\x,t)|\vv(\x,t)=-\nabla p(\x,t).
\eeq

Forchheimer's three-term law
\beq\label{F3}
a(\x,t)v+b(\x,t)|\vv(\x,t)|\vv(\x,t)+c(\x,t)|\vv(\x,t)|^2\vv(\x,t)=-\nabla p(\x,t). 
\eeq

Forchheimer's power law
\beq\label{FP}
a(\x,t) \vv(\x,t)+d(\x,t)|\vv(\x,t)|^{m-1}\vv(\x,t)=-\nabla p(\x,t).
\eeq
Here, the $a(\x,t)$, $b(\x,t)$, $c(\x,t)$, $d(\x,t)$ are positive functions, and $m\in(1,2)$ are derived from experiments for each case.
The above three Forchheimer equations can be combined and generalized to the following form:
\beq\label{gF}
g_F(|\vv(\x,t)|)\vv(\x,t)=-\nabla p(\x,t),
\eeq
where
\beq
g_F(y)=a_0(\x,t)+a_1(\x,t)y^{\alpha_1}+\dots+a_N(\x,t)y^{\alpha_N},
\eeq
with $N\ge 1$, $\alpha_0=0<\alpha_1<\ldots<\alpha_N$, $a_0(\x,t),a_N(\x,t)>0$, $a_1(\x,t),a_2(\x,t),\ldots,a_{N-1}(\x,t)\ge 0$.
The generalized Forchheimer equation \eqref{gF} was intensely used by the authors to model and study fast flows in the porous media, see in  \cite{ABHI1,HI1,HI2,HIKS1,HKP1,HK1,HK2,CHK1,CHK2}. The techniques developed in those papers will be essential in our approach and analysis below.

In previous works, each regime pre-Darcy, Darcy, or post-Darcy was studied separately, even though they  exist simultaneously in porous media. In particular cases, some models must consider multi-layer domains with each layer having a different regime of fluid flows, see for e.g. section 6.7.8 of \cite{StraughanBook}.
The goal of this section is to model all regimes together in the same domain. 
\newcommand{\Gfield}{\mathbf G}
We write  a general equation of motion for all cases \eqref{D}--\eqref{gF} as
\beq\label{vecform}
g(|\vv(\x,t)|)\vv(\x,t)=-\nabla p(\x,t),
\eeq
Different forms of $g(y)$ give different models, for example, 
\beqs
g(y)=k^{-1}(\x,t)y^{-\alpha},\
k^{-1}(\x,t),\
a(\x,t)+b(\x,t)y,\
a(\x,t)+b(\x,t)y+c(\x,t)y^2,\
a(\x,t)+d(\x,t)y^{m-1},\
g_F(y),
\eeqs
for equations \eqref{PD}, \eqref{D}, \eqref{F2}, \eqref{F3}, \eqref{FP}, \eqref{gF}, respectively.

We consider the following two  main  models. Below, $\indic_E$ denotes the characteristic (indicator) function of a set $E$.

\textbf{Model 1.} Function $g(y)$ is piecewise smooth on $(0, \infty)$. 
Based on \eqref{PD}, \eqref{D} and \eqref{gF} and their validity in different ranges  of Reynolds number, our first consideration is the following piecewise defined function
\beq \label{gdef}
g(y)=\bar{g}(y)\eqdef \underline{a} y^{-\alpha}\indic_{(0,y_1)}(y)+ \bar{a}\indic_{[y_1,y_2]}(s)+g_F(y)\indic_{(y_2,\infty)}(y)
\quad \text{ for } s>0,
\eeq
where $\alpha\in(0,1)$, and $y_2>y_1>0$ are fixed threshold values.

%
Note that $\bar g(y)$ is not differentiable at $y_1,y_2$.

\textbf{Model 2.} Function $g(y)$ is smooth on $(0, \infty)$.
Another generalization is  to use a smooth interpolation between pre-Darcy \eqref{PD} and generalized Forchheimer \eqref{gF}. Instead of \eqref{gdef}, we propose the following
\begin{equation}\label{gI}
g(y)=g_I(y)\eqdef a_{-1}y^{-\alpha}+a_0+a_1y^{\alpha_1}+\dots+a_Ny^{\alpha_N} 
\quad \text{for } y>0,
\end{equation}
where  $N\ge 1$, $\alpha\in (0,1)$, $\alpha_N>0$, 
\beq\label{aicond} 
a_{-1},a_0,a_N>0 \text{ and }a_i\ge 0\quad  \forall i=1,\ldots, N-1.
\eeq
The main advantage of $g_I$ over $\bar g$ is that it is smooth on $(0,\infty)$.
This  allows further mathematical analysis of the flows.
It also can be used as a framework for perspective interpretation of field data, i.e., matching the coefficients $a_i$ for $i=-1,0,1,\ldots,N$ to fit the data.   
The  Darcy--Forchheimer equation,  which is considered as a momentum equation, is studied in \cite{ABHI1,HI1,HI2} of the form 
\beq
-\nabla p(\x,t) =\sum_{i=-1}^N a_i(\x,t) |\vv(\x,t)|^{\alpha_i}\vv(\x,t). 
\eeq  
where  $N\ge 1$, $\alpha_{-1}=-\alpha, \alpha\in (0,1)$, $\alpha_N>0$,  
$a_{-1}(\x,t),a_N(\x,t)>0 \text{ and }a_i(\x,t)\ge 0\quad  \forall i=1,\ldots, N-1.
$

In order to take into account the presence of density in generalized Forchheimer equation, we modify \eqref{gF} using dimension analysis by Muskat \cite{Muskatbook} and Ward \cite{Ward64}. They proposed the following equation for both laminar and turbulent flows in porous media:
\beq\label{W}
-\nabla p(\x,t) =G\big(\vv^i \kappa^{\frac {i-3} 2} \rho^{i-1} \mu^{2-i}\big),
\eeq 
where  $G$ is a function of one variable, $\mu=\mu(\x,t)$ is the viscosity of the fluid, $\kappa=\kappa(\x)$ is  the permeability of the medium. 
  
  In particular, when $i = 1$, Ward \cite{Ward64} established the Darcy's law to match the experimental data
\beq\label{FW}
-\nabla p(\x,t)=\frac{\mu(\x,t)}{\kappa(\x)} \vv(\x,t),
\eeq
and when $i = 2$ for Forchhemer's law
 \beq\label{FW1}
-\nabla p(\x,t)=\frac{\mu(\x,t)}{\kappa(\x)} \vv(\x,t)+c_F\frac{\rho(\x,t)}{\sqrt {\kappa(\x)}}|\vv(\x,t)|\vv(\x,t),\quad \text{where }c_F>0.
\eeq

Combining  \eqref{vecform} with the suggestive form \eqref{W} for the dependence on $\rho$ and $\vv$, we propose the following equation 
 \beq\label{FM}
-\nabla p(\x,t)= \sum_{i=-1}^N a_i(\x,t) \rho^{\alpha_i} |\vv(\x,t)|^{\alpha_i} \vv(\x,t),
 \eeq
where $N\ge 1,\alpha_{-1}=-\alpha <\alpha_0=0<\alpha_1<\ldots<\alpha_N$ are fixed real numbers, the coefficients $a_{-1}(\x,t), \ldots, a_N(\x,t)$ are non-negative functions with 
$$ 0<a_* <a_{-1}(\x,t),a_0(\x,t) ,a_N(\x,t)<a^*<\infty, \quad 0\le a_i(\x,t) \le a^*<\infty,\,  i=2,\ldots, N-1.$$

Multiplying both sides of the equation \eqref{FM}  to $\rho$, we find that  
 \beq\label{eq1}
  \Big(\sum_{i=-1}^N a_i(\x,t) |\rho(\x,t) \vv(\x,t)|^{\alpha_i}\Big) \rho(\x,t) \vv(\x,t)   =-\rho(\x,t)\nabla p(\x,t).
 \eeq
 
 Denote the function $F:\Omega\times[0,T]\times\mathbb{R}^+\rightarrow\mathbb{R}^+$ a generalized polynomial with non-negative coefficients by
\beq\label{GPoly}
F(z)=a_{-1}(\x,t)z^{-\alpha}+  a_0(\x,t) + a_1(\x,t)z^{\alpha_1}+\cdots +a_N(\x,t)z^{\alpha_N},\quad z\ge 0. 
\eeq 

We denote the vector of powers by $\alp=(\alpha_{-1},\alpha_0, \alpha_1\ldots,\alpha_N)$
and vector coefficients by \\
 $\bfa(\x,t)=\left(a_{-1}(\x,t),a_0(\x,t), a_1(\x,t)\ldots,a_N(\x,t)\right).$ 
The class of functions $F(s)$ as in \eqref{GPoly} is denoted by P($N,\alp$).

The equation \eqref{eq1} can be rewritten as 
\beq\label{eq1a} 
F(|\rho(\x,t) \vv(\x,t)|)\rho(\x,t) \vv(\x,t) = -\rho(\x,t)\nabla p(\x,t).
\eeq

Above, we have introduced several models which can be used to interpret experimental and field data. We now use them to investigate the fluid flow's properties. They are used together with other basic equations of continuum mechanics which we recall here.

Continuity equation
\beqs
\phi(\x)\rho_t(\x,t)+\nabla\cdot(\rho(\x,t) \vv(\x,t))=0,
\eeqs
where $\phi\in(0,1)$ is the constant porosity.

Under isothermal condition the state equation relates the density $\rho$ with the pressure $p$ only, i.e., $\rho=\rho(p)$. Therefore, the equation of state, for slightly compressible fluids, is
\beqs
\frac{d\rho}{dp}=\frac{\rho}\kappa,
\eeqs
where $1/\kappa>0$ is small compressibility.

Hence
\beq \label{eq3}
\nabla \rho = \frac{1}{\kappa} \rho\nabla p, 
\quad \text{ or }\quad  \rho \nabla p=\kappa\nabla \rho.
\eeq
Combining \eqref{eq1a} and \eqref{eq3} implies that  
 \beq\label{ru}
 F(|\rho(\x,t) \vv(\x,t)|) \rho(\x,t) \vv(\x,t)  =-\kappa(\x)\nabla \rho(\x,t).
 \eeq

The continuity equation is
\beq\label{con-law}
\phi(\x)\rho_t(\x,t)+{\rm div }(\rho(\x,t) \vv(\x,t))=f(\x,t),
\eeq
where $\phi$ is the porosity, $f$ is external mass flow rate . 

By combining \eqref{ru} and \eqref{con-law} we have
\beqs
\begin{aligned}
&F(|\m(\x,t) |) \m(\x,t) = -\kappa(\x)\nabla \rho(\x,t), \\
&\phi(\x)\rho_t(\x,t)+ \diver{\m(\x,t)}=f(\x,t),
\end{aligned}
\eeqs
where $\m(\x,t)=\rho(\x,t) \vv(\x,t)$.

By rescaling the variable $\rho(\x,t)\to\kappa(\x)\rho(\x,t)$, $\phi(\x)\to \kappa^{-1}(\x)\phi(\x) $, we  obtain system of equations  
\beq\label{main-sys}
\begin{aligned}
&F(|\m(\x,t) |) \m(\x,t) = -\nabla \rho(\x,t), \\
&\phi(\x)\rho_t(\x,t)+ \diver{\m}(\x,t)=f(\x,t).
\end{aligned}
\eeq  
 The outline of the paper is as follows.  
 We introduce the notations and the relevant results in section \S \ref{Intrsec}.
 In section \S \ref{MFEM}, we defined a numerical approximation using mixed finite element approximations 
 to the initial boundary value problem (IVBP)~\eqref{MainProb}. 
 In section \S \ref{Bsec}, we establish many estimates of the energy type norms for the approximate solution $(\m_h, \brho_h)$ to the IVBP problem~\eqref{semidiscreteform} in Lebesgue norms in terms of the boundary data and the initial data.
In section \S \ref{dependence-sec}, we focus on proving the continuous dependence of the solution on the coefficients of function $g$. 
It is proved in Theorem~\ref{DepCoeff} that the difference between the two solutions, which corresponds to two different coefficient vectors $\bfa_1$ and  $\bfa_2$ is estimated in terms of $|\bfa_1-\bfa_2|$, see in \eqref{ssc2}.
In section \S \ref{err-sec}, we study in Theorem~\ref{err-thm} the convergence and in Theorem~\ref{par-thm} the dependence on coefficients of general polynomial $g$ of the approximated solution to the problem~\eqref{semidiscreteform}. Furthermore, we can specify the convergent rate.          
In section \S \ref{err-ful-sec}, we study the fully discrete version of problem~\eqref{semidiscreteform}. In Lemma~\ref{stab-appr}, the stability of the approximated solution is proved. Theorems~\ref{Err-ful}~and~\ref{Dep-ful} are for studying the error estimates and the continuous dependence on parameters of the numerical solution.   
In section \S \ref{Num-result}, the numerical experiments in the two-dimensions using the standard finite elements $\mathbb P_1$ are presented regarding the convergence rates and the dependence of the solution on the physical parameters.
\section{Notations and preliminary results}\label{Intrsec}

Through out this paper, we assume that $\Omega$ is an open, bounded subset of $\mathbb{R}^d$, with $d=2,3,\ldots$, and has $C^1$-boundary $\partial \Omega$. For $s\in [0,\infty)$, we denote $L^{s}(\Omega)$ be the set of s-integrable functions on $\Omega$ and $( L^{s} (\Omega))^d$ the space of $d$-dimensional vectors which have all components in $L^{s}(\Omega)$.  We denote $(\cdot, \cdot)$ the inner product in either $L^{s}(\Omega)$ or $(L^{s}(\Omega))^d$ that is
$
( \xi,\eta )=\int_\Omega \xi\eta dx$  or $(\boldsymbol{\xi},\boldsymbol \eta )=\int_\Omega \boldsymbol{\xi}\cdot \boldsymbol{\eta} dx 
$ 
and $ \norm{v}_{L^s(\Omega)}=\left(\int_{\Omega} |v(x)|^s dx\right)^{1/s}$ for standard  Lebesgue norm of the measurable function. 
 For $m\ge 0, s\in [0,\infty]$, we denote the Sobolev spaces by $W^{m,s}(\Omega)=\{v \in L^s(\Omega): D^{\alpha} v \in L^s(\Omega), |\alpha|\le m \}$ and the norm of $W^{m,s}(\Omega)$ by $
 \norm{v}_{W^{m,s}(\Omega)} = \left(\sum_{|\alpha|\le m}\int_{\Omega} |D^\alpha v|^s dx\right)^{1/s},$  and  $\norm{v}_{W^{m,\infty}(\Omega)} = \sum_{|\alpha|\le m}\esssup_{\Omega} |D^\alpha v|.$ Let $I=[0,T]$, we define $L^s(I,X)$ to be the space of all measurable functions $v : I\to X$ with the norm $\norm{v}_{L^s(I,X)}=\left(\int_0^T \norm{v(t)}_X^s dt\right)^{1/s}$,  and $L^\infty(I;X)$ to be the space of all measurable functions $v: I \to X$ such that $v: t\to \norm{v(t)}_X$ is essentially bounded on $I$ with the norm $\norm{v}_{L^\infty(I,X)}=\esssup_{t\in I}\norm{v(t)}_X$. 
We use short hand notations, 
$
\norm{\rho(t)} = \norm{ \rho(\cdot, t)}_{L^2(\Omega)}, \forall t\ge 0$  and $\rho_0(\cdot) =  \rho(\cdot,0).$

Our estimates make use of coefficient-weighted norms. For some strictly positive,
bounded function, we denote the weighted $L^2$-norm by $\norm{f}_\omega$ by 
$
\norm{f}_{\omega}^2 =\int_\Omega \omega |f|^2 dx,  
$
and if $0<\omega_* \le \omega (x) \le \omega^*$ throughout $\Omega$, we have the equivalent 
$
\sqrt{\omega_*} \norm{f}_{L^2(\Omega)} \le \norm{f}_{\omega} \le \sqrt{\omega^*}\norm{f}_{L^2(\Omega)}. 
$
Our calculations frequently use the following exponents
\beq\label{a-const }
   s= \alpha_N+2,\quad  s^*=\frac{s}{s-1}.
  \eeq
The argument $C$ will represent for positive generic constants and their values  depend on exponents, coefficients of polynomial  $F$,  the spatial dimension $d$ and domain $\Omega$, independent of the initial and boundary data and time step. These constants may be different place by place. 


We will frequently use the following basic inequalities.
By Young's inequality, we have
\begin{align}
\label{bi3}
x^\beta&\le x^{\gamma_1}+x^{\gamma_2}\quad\text{for all }x>0,\ \gamma_1\le \beta \le \gamma_2,\\
\label{bi2}
x^\beta&\le 1+x^\gamma\quad\text{for all }x\ge 0, \ \gamma\ge\beta > 0.
\end{align}
For any $r\ge1$, $x_1,x_2,\ldots,x_k\ge 0$, 
\beq\label{bi0}
x_1^r+x_2^r+\cdots+x_k^r \le (x_1+x_2+\cdots+x_k)^r\le k^{r-1}(x_1^r+x_2^r+\cdots+x_k^r). 
\eeq

\begin{lemma} For all $w>0$,

(i) 
 \beq\label{dervF}
-\alpha F(w)  \le wF'(w)\le \alpha_N F(w).
 \eeq
 
 (ii) 
\beq\label{OrdF}
 a_* (w^{-\alpha}+1+w^{\alpha_N}) \le F(w)\le N a^* (w^{-\alpha}+w^{\alpha_N}). 
\eeq

\end{lemma}
\begin{proof}

(i)
Since $-\alpha<\alpha_0\le\alpha_1\le \ldots\le\alpha_N$, 
\beqs
  -\alpha F(w) \le  w F'(w) = -\alpha a_{-1}w^{-\alpha}+\alpha_0a_0w^{\alpha_0}+ \alpha_1a_1w^{\alpha_1} +\cdots+ \alpha_{N}a_N w^{\alpha_N} \le \alpha_N F(w).
 \eeqs 
 Thus, the estimate \eqref{dervF} follows. 

(ii) It is easy to see that
\beq\label{upineq}
F(w)\le \max_{i=-1,\ldots,N}a_i\cdot (w^{-\alpha}+w^{\alpha_0}+\cdots + w^{\alpha_i}) 
\le N a^* (w^{-\alpha}+w^{\alpha_N}).
\eeq
On the other hand, 
\beq\label{lowineq}
F(w)\ge a_{-1}w^{-\alpha}+ a_0 +a_Nw^{\alpha_N} 
\ge a_* (w^{-\alpha}+1+ w^{\alpha_N}).
\eeq
Combining \eqref{upineq} and \eqref{lowineq} we obtain the inequality in \eqref{OrdF}. 
\end{proof}
\begin{lemma}
(i)  Assume  $-1<p <0$, then for all $\x,\y\in \R^d$,    
\beq\label{cont1}
||\x|^p \x -|\y|^p \y |\le 2|\x-\y|^{1+p}.
\eeq
(ii)  Assume $p>0$, then for all $\x,\y\in \R^d$,  there is $C>0$ such that  (see in \cite{M2AN75,SJ78}). 
\beq\label{cont2}
||\x|^p \x-|\y|^p \y | \le C(|\x|+ |\y|)^p |\x-\y|. 
\eeq
\end{lemma} 
\begin{proof}
 
 (i) Let $\gamma(\tau)=\tau  \x+ (1-\tau)\y, \tau\in [0,1]$ and $h(\tau) =|\gamma(\tau)|^p\gamma(\tau)$. Then 
\beq\label{asw}
\begin{split} 
\left| |\x|^p \x -|\y|^p \y  \right|&= |h(1)-h(0)| \\
&=\left|\int_0^1 h'(\tau) d\tau\right|
=\left|\int_0^1 \Big[|\gamma(\tau)|^p (\x-\y) + p|\gamma(\tau)|^{p-1} \frac{\gamma(\tau)\cdot(\x-\y)}{|\gamma(\tau)|}\gamma(\tau)\Big]   d\tau \right|\\
&\le (1+p)|\x-\y|\int_0^1 |\gamma(\tau)|^p   d\tau.
\end{split} 
\eeq
We claim  
\beq\label{estgamma}
 \int_0^1 |\gamma(\tau)|^p   d\tau\le \frac 2{p+1}|\x-\y|^p. 
\eeq
 The inequality \eqref{cont1} follows by substituting \eqref{estgamma} into \eqref{asw}. 
 
{ \bf Proof of claim \eqref{estgamma}}

 \textbullet~~ If $|\x|\ge |\x-\y|$. Since $p<0$, 
\begin{align*}
|\tau \x+(1-\tau)\y|^p\le ||\x| -(1-\tau)|\x-\y||^p
\le | (1-\tau)|\x-\y| - |\x-\y| |^p = \tau^p|\x-\y|^p.   
\end{align*}
This shows that 
\beqs
\int_0^1 |\gamma(\tau)|^p d\tau \le |\x-\y|^p \int_0^1 \tau^p d\tau \le  \frac{2}{p+1} |\x-\y|^p. 
\eeqs 

\textbullet~~ If $|\x|< |\x-\y|$. Let $\tau_*\in(0,1)$ be defined by $(1-\tau_*)|\x-\y| =|\x|$
\begin{align*}
\int_0^1 |\gamma(\tau)|^p d\tau &\le \int_0^1 | |\x| -(1-\tau)|\x-\y| |^p d\tau\\
&=|\x-\y|^p \int_0^1 | \tau-\tau_*|^p d\tau
=\frac 1{p+1} |\x-\y|^p( \tau_*^{p+1}+ (1-\tau_*)^{p+1} ) \\
&\le \frac 2{p+1} |\x-\y|^p.
\end{align*}
(ii) For the proof see Lemma  5.3 in \cite{M2AN75,SJ78}. 
\end{proof}

For convenience, we use the following notations: let $\y=(y_1, y_2, \ldots, y_d)$ and $\y'=(y'_1, y'_2, \ldots, y'_d)$ be two arbitrary vectors of the same length, including possible length $1$.
\begin{lemma} Let $F(\bfa,\cdot)$ and $F(\bfa',\cdot)$ belong to the class P($N,\alp$). Then for any  $\y,\y'\in \mathbb R^d$, one has 

(i)
\beq\label{Umono}
|F(\bfa,|\y|) \y-F(\bfa',|\y'|)\y'|\le  C_1(1+|\y|+|\y'|)^{s-2+\alpha}|\y-\y'|^{1-\alpha}
+ C_2 (1+|\y|+|\y'|)^{s-1}|\bfa -\bfa' |. 
\eeq

(ii)
\beq\label{quasimonotone}
(F(\bfa,|\y|) \y-F(\bfa',|\y'|)\y')\cdot (\y-\y')
\ge C_3 |\y'-\y|^s
 - C_2(1+|\y|+|\y'|)^{s-1} |\y-\y'| \,|\bfa -\bfa' |, 
\eeq
where   
$
C_1=2(s-1)/(1-\alpha),C_2=3N ,C_3 =a_*(1-\alpha)/(2^{s-1}(s-1)).
$

In particular, in the case $\bfa = \bfa'$ then \eqref{Umono} and \eqref{quasimonotone} become

(iii)
\beq\label{Lipchitz}
   |F(\bfa,|\y|) \y-F(\bfa,|\y'|)\y'|\le  C_1 (1+|\y|+|\y'|)^{s-2+\alpha}|\y-\y'|^{1-\alpha}.
\eeq  
 
(iv)
\beq\label{monotone0}
   (F(\bfa,|\y|) \y-F(\bfa,|\y'|)\y')\cdot (\y-\y') \ge C_3 |\y'-\y|^{s}.
\eeq 
\end{lemma}
\begin{proof} 
(i). Let  $\bfa$, $\bfa'\in S$ and  $\y,\y' , \kk \in \mathbb R^d$.
We defined,  
\beqs
\gamma(\tau)=\tau \y +(1-\tau)\y', \quad
\mathbf{b}(\tau)=\tau\bfa+(1-\tau)\bfa' \quad \forall \tau\in [0,1],\\
\eeqs

{\it Case 1}: The origin does not belong to the segment connecting  $\y'$ and  $\y$. Define 
 $$z(\tau)=F(\mathbf{b}(\tau),|\gamma(\tau)|)\, \gamma(\tau)\cdot \kk .$$
 By the Mean Value Theorem,  we have 
\beqs
I
\eqdef [F(\bfa,|\y|) \y-F(\bfa',|\y'|)\y']\cdot \kk
=z(1)-z(0)=\int_0^1 z'(\tau) d\tau.
\eeqs
Elementary calculations give
\beq\label{Ifm}
I=\int_0^1 f_1(\tau) d\tau + \int_0^1 f_2(\tau) d\tau \eqdef I_1+I_2,
\eeq
where
\begin{align*}
f_1(\tau)&=F(\mathbf{b}(\tau),|\gamma(\tau)|)(\y-\y')\cdot \kk
 +F'(\mathbf{b}(\tau),|\gamma(\tau)|)\frac{\gamma(\tau)\cdot(\y-\y')}{|\gamma(\tau)|} \gamma(\tau)\cdot \kk,\\
f_2(\tau)&= F_{\bfa}( \mathbf{b}(\tau),|\gamma(\tau)|)(\bfa-\bfa')\gamma(\tau)\cdot \kk .
\end{align*}
 $\bullet$ \emph{Estimation of $I_1$.}  Using triangle inequality, \eqref{dervF} and \eqref{OrdF} we find that 
\beqs
\begin{split}
|f_1(\tau)| &\le \left |F(\mathbf{b}(\tau),|\gamma(\tau)|)(\y-\y')
 +F'(\mathbf{b}(\tau),|\gamma(\tau)|)\frac{\gamma(\tau)\cdot(\y-\y')}{|\gamma(\tau)|} \gamma(\tau)  \right|\,  |\kk| \\
 &\le (1+\alpha_N)F(|\mathbf{b}(\tau),\gamma(\tau)|) |\y-\y'|\, |\kk|\\
 &\le (1+\alpha_N)(|\gamma(\tau)|^{-\alpha}+|\gamma(\tau)|^{\alpha_N})  |\y-\y'|\, |\kk|\\ 
 &\le (1+\alpha_N)|\gamma(\tau)|^{-\alpha}(1+|\y|+|\y'|)^{\alpha_N+\alpha}  |\y-\y'|\, |\kk|.
\end{split}
\eeqs
Thus 
\beqs
I_1 \le \int_0^1 |f_1(\tau)| d\tau \le  (1+\alpha_N)(1+|\y|+|\y'|)^{\alpha_N+\alpha}|\y-\y'|\, |\kk|\, \int_0^1 |\gamma(\tau)|^{-\alpha} d\tau.
\eeqs
By using \eqref{estgamma} yields    
\beq\label{UI1}
I_1 \le \frac{2(1+\alpha_N)}{1-\alpha} (1+|\y|+|\y'|)^{\alpha_N+\alpha}|\y-\y'|^{1-\alpha}\, |\kk|.
\eeq
$\bullet$ \emph{Estimation of $I_2$.} We find the partial derivative of $F(\bfa, |\xi|)$ in $\bfa$.
For $i=-1,0,\ldots,N$,  taking the partial derivative in $a_i$, we find that
\beq\label{Ka}
 \sum_{i=-1}^N F_{a_i}(\bfa,|\xi|)
 = |\xi|^{-\alpha} + |\xi|^{\alpha_0}+\cdots + |\xi|^{\alpha_N}
 \le N(|\xi|^{-\alpha}+1+|\xi|^{\alpha_N})\quad \text{ by \eqref{bi2}}.
\eeq
Using the estimate \eqref{OrdF}, we bound 
\beq\label{ih2}
\begin{aligned}
| f_2(\tau)|
&\le |F_{\bfa}(\mathbf{b}(\tau),|\gamma(\tau)|)|\ |\bfa -\bfa'|\  | \gamma(\tau) |\ |\kk|\\
&\le N(|\gamma(\tau)|^{-\alpha}+1+|\gamma(\tau)|^{\alpha_N})\ |\bfa -\bfa'|\  | \gamma(\tau) |\ |\kk|\\
&\le 3N(1+|\gamma(\tau)|^{1+\alpha_N})  \  |\bfa -\bfa'|  \ |\kk|.
\end{aligned}
\eeq
Using the fact $|\gamma(\tau)|\le |\y|+ |\y'|$,  \eqref{ih2} yields
\beq\label{Uh2}
| f_2(\tau)| \le 3N(1+|\y|+|\y'|)^{\alpha_N+1}  \  |\bfa -\bfa'|  \ |\kk|,
\eeq
and consequently,
\beq \label{UI2}
I_2\le  \int_0^1 | f_2(\tau)|d\tau 
\le  3N(1+|\y|+|\y'|)^{\alpha_N+1}|\bfa -\bfa' | \ |\kk|. 
\eeq

Then \eqref{Umono} follows by combining \eqref{Ifm}, \eqref{UI1} and \eqref{UI2}.

{\it Case 2}: The origin belongs to the segment connect $\y', \y$. We replace $\y'$ by $\y'+\y_\varep$ so that $0\notin [\y'+\y_\varep, \y]$ and $y_\varep \to 0$ as $\varep \to 0$. Apply the above inequality for $\y$ and $\y'+\y_\varep$, then let $\varep\to 0$.

{(ii)} If $\kk = \y-\y'$ then 
\beq\label{ks}
\begin{aligned}
f_1(\tau) 
&\ge F(\mathbf{b}(\tau),|\gamma(\tau)|)|\y-\y'|^2 - \alpha \frac{F(\mathbf{b}(\tau),|\gamma(\tau)|)|\, |\gamma(\tau)\cdot(\y-\y')|^2}{|\gamma(\tau)|^2} \\
&\ge (1-\alpha) F(\mathbf{b}(\tau),|\gamma(\tau)|)|\y-\y'|^2.
\end{aligned}
\eeq

It follows \eqref{ks} and \eqref{estgamma} that 
\beq\label{int-h1}
\begin{split}
I_1
  &\ge  (1-\alpha)|\y'-\y|^2 \int_0^1 F(\mathbf{b}(\tau),|\gamma(\tau)|) d\tau\\
  &\ge a_*(1-\alpha )|\y'-\y|^2\int_0^1 (|\gamma(\tau)|^{-\alpha} +1+|\gamma(\tau)|^{\alpha_N})  d\tau\\
  & > a_*(1-\alpha )|\y'-\y|^2\int_0^1 |\gamma(\tau)|^{\alpha_N}  d\tau.
\end{split}
\eeq
It is proved in Lemma 2.4 of \cite{CHIK1}  that  , 
\beq
\int_0^1 |\gamma(\tau)|^{\alpha_N} ds \ge \frac{|\y'-\y|^{\alpha_N}}{2^{\alpha_N+1}(\alpha_N+1) }. 
\eeq
Substituting this into \eqref{int-h1} gives 
\beq\label{mono-est}
I_1 \ge \frac{a_*(1-\alpha ) }{2^{\alpha_N+1}(\alpha_N+1) } |\y'-\y|^{\alpha_N+2}.
\eeq
On the other hand from \eqref{ih2}, we see that
\beq \label{I2geq}
I_2\ge - \int_0^1 | f_2(\tau)| d\tau 
\ge - 3N (1+|\y|+|\y'|)^{\alpha_N+1} |\bfa -\bfa' | \ |\y-\y'|.
\eeq
Thus, we obtain \eqref{quasimonotone} by combining \eqref{Ifm}, \eqref{int-h1} and \eqref{I2geq}.
\end{proof}

We recall a discrete version of Gronwall Lemma in backward difference form, which is useful later. 
\begin{lemma}[Lemma 2.4 in~\cite{K17}]
\label{DGronwall}
	Assume the nonnegative sequences $\{a_n\}_{n=0}^\infty$, $\{b_n\}_{n=0}^\infty$ ,$\{g_n\}_{n=0}^\infty$  satisfying 
	\beqs
	\frac{a_n-a_{n-1}}{\Dt} - a_n  + b_n \le g_n, \quad n=1,2,3\ldots
	\eeqs
	then for sufficiently small $\Dt$, 
	\beq\label{Gineq}
	a_n   +\Dt \sum_{i=1}^n b_i  \le e^{\frac{n\Dt}{1-\Dt}}     \Big(  a_0 +\Dt\sum_{i=1}^{n} g_i \Big).
	\eeq
\end{lemma}
%
%

\section{A mixed finite element method  and mixed finite element approximation}\label{MFEM}
In this section we will represent the weak formulation and mixed element approximation of the mixed regime equation. We consider the following model problem   
\beq\label{TProb}
\begin{aligned}
F(|\m(\x,t) |) \m(\x,t) = -\nabla    \rho(\x,t) & & \quad (\x,t)\in\Omega\times I , \\
\phi(\x)  \rho_t(\x,t) + \nablad \m(\x,t)=f(\x,t) & & \quad(\x,t)\in\Omega\times I,\\
  \rho(\x,t) =\psi(\x,t) & & \quad (\x,t)\in\partial\Omega \times I ,\\
\rho(\x,0)=\rho_0(\x) & & \quad \x\in\Omega.
\end{aligned}
\eeq
We make the following assumptions on the data and coefficients 
 \begin{itemize}
 \item [(H1)] The coefficients $a_{-1}, \ldots, a_N\in W^{1,\infty}(I,L^\infty(\Omega)) $ satisfy  
\begin{align*} 
&0<a_* <a_{-1}(\x,t),a_0(\x,t), a_N(\x,t)<a^*<\infty,  0\le a_i(\x,t)\le a^*<\infty,\,  i=1,\ldots, N-1 \\
 &|a_{it}(\x,t)|< b^* <\infty , i=-1,\ldots, N
 \end{align*}
for almost every  $(\x,t)\in \bar\Omega\times I$. 

 \item[(H2)] $0<\phi_*\le \phi(\x)\le \phi^*<\infty;
 f \in L^\infty (I,L^2(\Omega) );
\psi\in C_{\x,t}^{1,2} (  \bar\Omega \times I  ).$  
 \end{itemize}
Dealing with the boundary condition, let $\Psi(\x,t)$ be an extension of $\psi$ from $\partial \Omega\times I$ to $ \bar \Omega\times I$ (see e.g \cite{HI1,JS1995}).
Let $\bar \rho =\rho -\Psi$. Then
\beq\label{MainProb}
\begin{aligned}
	F(|\m(\x,t) |) \m(\x,t) = -\nabla\brho(\x,t) -\nabla \Psi(\x,t) & & \quad (\x,t)\in\Omega\times I , \\
	\phi(\x) \brho_t(\x,t) + \nablad \m(\x,t)=f(\x,t) -\phi(\x)\Psi_t(\x,t)  & & \quad(\x,t)\in\Omega\times I,\\
	\brho(\x,t) =0 & & \quad (\x,t)\in\partial\Omega \times I ,\\
	\brho(\x,0)=\rho_0(\x)-\Psi(\x,0) & & \quad \x\in\Omega.
\end{aligned}
\eeq
Define $Q=L^2(\Omega)$, and 
 $
 V=\left\{ \vv\in (L^s(\Omega))^d , \nabla\cdot\vv \in L^2(\Omega)\right\} 
 $
 and equip it with the norm
 \beq \label{Greenthm}
 \norm{\vv}_{V} = \norm{\vv}_{L^s(\Omega)} + \norm{\nabla\cdot \vv}_{L^2(\Omega)}.
  \eeq
The mixed formulation of \eqref{TProb} is defined as follows:
 Find $(\m,\bar\rho): I\to  V\times Q$ such that  
\beq\label{WeakProb}
\begin{aligned}
 \intd{ F(|\m |) \m, \vv} -(  \bar \rho, \nabla\cdot \vv)= -\intd{ \nabla\Psi,\vv }, \quad \text{ for all } \vv\in V,\\
\intd{\phi \bar \rho_t, q} + \intd{\nablad \m , q} =\intd{ f , q} -\intd{\phi\Psi_t, q},  \quad\text{ for all } q\in Q.
\end{aligned}
\eeq

 Let  $\{\mathcal T_h\}_h$ be a family of globally quasi-uniform triangulation of $\Omega$ with $\max_{\tau\in \mathcal T_h} \rm{ diam}~\tau \le h$. Let $k\ge 0$  be an integer we define 
\begin{align*}
Q_h &=\{\rho_h\in C^0(\overline \Omega), \forall \tau\in \mathcal T_h, \rho_h|_\tau \in P_{k}(\tau)\}, \\
V_h&=\{\m_h\in C^0(\overline \Omega)^d, \forall \tau\in \mathcal T_h, \m_h|_\tau \in P_{k}(\tau)\},
\end{align*}
  with $P_k(\tau)$ being the space of polynomial of degree at most $k$ on the element $\tau$.    Let $V_h \times Q_h$ be the mixed element spaces approximating the space  $V\times Q$.     
  
 For momentum, let $\Pi: V\to V_h$ be the Raviart-Thomas projection \cite{RT06}, which satisfies
 \beqs
   \intd{\nabla\cdot( \Pi \m -\m), q} =0, \quad \text {for all } \m\in V, q\in Q_h. 
 \eeqs     
 For density, we use the standard $L^2$-projection operator, see in \cite{Ciarlet78},  
$\pi: Q \to Q_h$,  satisfying
\begin{align*}
( \pi \rho -\rho, q ) = 0,& \quad \text {for all } \rho\in Q,  q\in Q_h,\\
( \pi \rho -\rho, \nablad \m_h ) = 0,& \quad\text {for all } \m_h\in V_h, \rho\in Q. 
\end{align*}

This projection has well-known approximation properties, e.g.~\cite{BF91,JT81,BPS02}.
\begin{align}
  &\label{DProj}  \norm{\Pi \m}_{L^q(\Omega)} \le C\left(\norm{\m}_{L^p(\Omega)}+h\norm{\nablad\m}_{L^2(\Omega)}\right),  \quad\forall \m\in V\cap (W^{1,q}(\Omega))^d.\\
\label{divProj}   & \norm{\Pi \m -\m }_{L^q(\Omega)} \le C h^p\norm{\m}_{W^{p,q}(\Omega)}, \quad  1/q <p\le k+1, \forall \m\in V\cap (W^{p,q}(\Omega))^d.\\
   & \norm{\nabla\cdot(\Pi \m -\m) }_{L^2(\Omega)}\le C h^p \norm{\nabla\cdot\m}_{W^{p,2}(\Omega) }, \quad  0  \le p\le k+1, \forall \m\in V\cap (H^p(\mathrm{div}, \Omega))^d. \\
&\label{LProj}\norm{\pi \rho}_{L^2(\Omega)}\le C \norm{\rho}_{L^2(\Omega)},  \quad \forall\rho\in L^2(\Omega).\\
\label{L2Proj} &\norm{\pi \rho - \rho }_{L^q(\Omega)} \leq C h^p \norm{\rho}_{W^{p,q}(\Omega)}, \quad  0\le p \le k+1,  q\in[1,\infty],\forall \rho\in W^{p,q}(\Omega).    
\end{align}
The two projections $\pi$ and $\Pi$ preserve the commuting property $\mathrm{div }\circ \Pi =\pi\circ \mathrm{ div}:V\to Q_h$. 

   The semidiscrete formulation of~\eqref{WeakProb} can read as follows: 
 Find a pair $(\m_h,\bar\rho_h): I \to V_h \times Q_h $ such that 
\beq\label{semidiscreteform}
\begin{aligned}
\intd{ F(|\m_h |) \m_h, \vv} -(  \brho_h, \nabla\cdot \vv)= -\intd{ \nabla\Psi,\vv } \quad \text{ for all } \vv\in V_h,\\
\intd{\phi \brho_{ht}, q} + \intd{\nablad \m_h , q} =\intd{ f , q} -\intd{\phi\Psi_t, q}  \quad\text{ for all } q\in Q_h.
\end{aligned}
\eeq
with initial data $\brho_{h0}=\pi\brho(\x,0)$.
\section{Stability}\label{Bsec}

We study the  equations \eqref{WeakProb}, and \eqref{semidiscreteform} for the density with fixed functions $F(s)$ in \eqref{eq1} and \eqref{GPoly}. 
Therefore, the exponents $\alpha_i$ and coefficients $a_i$ are all fixed, and so are the functions $F(\xi)$  in \eqref{GPoly}.  
With the properties \eqref{dervF}, \eqref{OrdF}, \eqref{Lipchitz}, the monotonicity \eqref{monotone0}, and by
classical theory of monotone operators \cite{MR0259693,s97,z90}, the authors in \cite {PG16, K1, K2} proved the global existence and uniqueness of the weak solution of the problem \eqref{semidiscreteform}.  

\subsection{A priori estimates for the solutions of the semi-discrete problems}
For  the { \it priori } estimates, we consider weak solutions with enough regularities in both  $\x$ and $t$ variables, but not necessarily classical, so
that our calculations can be applied. 

We will focus estimates  for $\bar\rho_h$. The estimates for $\rho_h$ can be obtained by simply using triangle inequality $\norm{\rho_h(t)}\le \norm{\bar\rho_h(t)}+\norm{\Psi(t)}$. Also our results are stated in term of $\Psi(x,t)$. These can be written in term of $\psi(x,t)$ by using specific extension. e.g. harmonic extension in \cite{HI1}.

\begin{lemma}
Let $(\m_h, \brho_h)$ be a solution to the problem \eqref{semidiscreteform}. There exists a positive constant $C$ independence of $h$ such that
\beq\label{stab1}
\norm{\brho_h }_{L^\infty(I,L^2(\Omega))}^2+\norm{\m_h}_{L^s(I,L^s(\Omega))}^s \le \norm{\brho(0)}_{L^2(\Omega)}^2 + C\mathscr{A}, 
\eeq
where
\beq\label{Adef}
\mathscr{A}= \norm{f}_{L^2(I,L^2(\Omega)) }^{2}+\norm{\Psi_t}_{L^2(I,L^2(\Omega))}^2+\norm{\nabla\Psi}_{L^{s^*}(I,L^{s^*}(\Omega))}^{s^*}.
\eeq 
\end{lemma}
\begin{proof}
Choosing $(\vv,q) = (\m_h, \brho_h)$ in \eqref{semidiscreteform} and adding the resultant equations yield
  \beq\label{deq}
   \frac 1 2 \ddt \norm{\brho_h}_{\phi}^2+\intd{ F(|\m_h|) \m_h, \m_h} =(f,\brho_h) - \intd{\phi\Psi_t,\brho_h}-\intd{ \nabla\Psi,\m_h }.
  \eeq
    By the monotonicity of the function $F(\cdot)$  in \eqref{monotone0}, the second term of \eqref{deq} is  bounded from below by 
  \beq\label{fest}
  \intd{ F(|\m_h|) \m_h, \m_h} \ge C_3\norm{\m}_{L^s}^s. 
  \eeq
 We bound the right hand side of \eqref{deq} by using Young's inequality to obtain  
  \beq\label{sest}
  (f,\brho_h)- \intd{\phi\Psi_t,\brho_h}-\intd{ \nabla\Psi,\m_h } \le  \norm{f}_{\phi^{-1}}^2 + \norm{\Psi_t}_{\phi}^2 + \frac {1}{2}\norm{\brho_h}_{\phi}^2+ \frac{C_3}{2} \norm{\m_h}_{L^s}^s + \frac{(s C_3/2)^{-s^*/s}}{s^*}  \norm{\nabla\Psi}_{L^{s^*}}^{s^*}.
  \eeq
Combining \eqref{deq}, \eqref{fest} and \eqref{sest} yields  
  \beqs
  \ddt \norm{\brho_h}_{\phi}^2+C_3 \norm{\m_h}_{L^s}^s 
  \le \norm{\brho_h}_{\phi}^2 + C\left(\norm{f}_{\phi^{-1}}^2 + \norm{\Psi_t}_{\phi}^2+\norm{\nabla\Psi}_{L^{s^*}}^{s^*}\right).
   \eeqs
  By Gronwall's inequality  and using the boundedness of the function $\phi$, we find that  
  \beqs
   \norm{\brho_h }_{L^\infty(I,L^2)}^2+\norm{\m_h}_{L^s(I,L^s)}^s 
  \le \norm{\brho_h(0)}_{L^2}^2 +C\left( \norm{f}_{L^2(I,L^2) }^{2}+ \norm{\Psi_t}_{L^2(I,L^2)}^2+\norm{\nabla\Psi}_{L^{s^*}(I,L^{s^*})}^{s^*}\right).
  \eeqs
Note that $\norm{\brho_h(0)}=\norm{\pi\brho(0)}\le \norm{\brho_0}$. Thus the inequality \eqref{stab1} holds.  
\end{proof}
\begin{theorem}
There is unique a solution of the problem \eqref{semidiscreteform} satisfying \eqref{stab1}.
\end{theorem}
\begin{proof}
 The equation \eqref{semidiscreteform} can be interpreted as the finite system of ordinary differential equations in the coefficients of $(\m_h, \bar\rho_h)$ with respect to basis of $V_h\times Q_h$. The stability estimates \eqref{stab1} suffice to establish the local existence of $(\m_h(t), \bar\rho_h(t))$ for all $t\in I.$ The proof of this statement is essential identical with that of \cite{EJP05, KIM199675} for generalized Forchheimer flows. We will omit here.      

Assume that $(\m_{h}^{(i)}, \brho_{h}^{(i)})$, $i=1,2$ are two solutions of \eqref{semidiscreteform}.  Let 
  $  \m_h= \m_h^{(1)} - \m_h^{(2)},  \brho_h = \brho_h^{(1)} - \brho_h^{(2)}.  $  Then
  \beq\label{DiffErr}
\begin{split}
\intd{ F(|\m_h^{(1)}|) \m_h^{(1)}- F(|\m_h^{(2)}|) \m_h^{(2)}, \vv} -\intd{\brho_h, \nabla\cdot \vv}  =0, \quad \forall\vv\in V_h,\\
\intd{ \phi \brho_{ht},q} + \intd{ \nabla\cdot \m_h, q}=0, \quad \forall q\in Q_h.
\end{split}
\eeq 
 It is easy to see that with $(\vv,q)=(\m_h,\brho_h)$ in \eqref{DiffErr}
    \beqs
\intd{ \phi \brho_{ht},\brho_h} +\intd{ F(|\m_h^{(1)}|) \m_h^{(1)}- F(|\m_h^{(2)}|) \m_h^{(2)}, \m_h} =0. 
\eeqs
Thanks to the monotonicity \eqref{monotone0}, we see that
\beq\label{a}
\frac 1 2\ddt \norm{\brho_h}_{\phi}^2 +C_3 \norm{\m_h}_{L^s}^s  \le \intd{ \phi \brho_{ht},\brho_h} +\intd{ F(|\m_h^{(1)}|) \m_h^{(1)}- F(|\m_h^{(2)}|) \m_h^{(2)}, \m_h} =0. 
\eeq 
This implies  
$
\norm{\brho_h}_{\phi}^2 +\norm{\m_h}^s_{L^s(I,L^s)} \le C\norm{\brho_h(0)}_{\phi}^2 =0. 
$
Hence $\brho_h=0$ and $\m_h=0$ a.e.    
\end{proof}
\begin{lemma} Let $(\m_h, \brho_h)$ be a solution to the problem \eqref{semidiscreteform}. Then, there exists a positive constant $C$ independence of $h$ such that
	\beq\label{stab2}
		\norm{\brho_{ht}}_{L^\infty(I,L^2(\Omega))}^2+\norm{\m_{ht}}_{L^2(I,L^2(\Omega))}^2 +\norm{\nablad \m_h}_{L^\infty(I,L^2(\Omega))}^2\le C (\mathscr{B}_0 +\mathscr{B}),
	 \eeq
	 where 
	 \beq\label{Bdef}
	 \mathscr{B}_0 =\norm{ \brho_0}_{L^2(\Omega)}^2+\norm{ \brho_t(0)}_{L^2(\Omega)}^2,\text { and }
	 \mathscr{B}=\mathscr{A} +\norm{f_t}_{L^2(I,L^2(\Omega))}^2+\norm{\Psi_{tt}}_{L^2(I,L^2(\Omega))}^2+\norm{\nabla\Psi_t}_{L^2(I,L^2(\Omega))}^2.
	 \eeq
\end{lemma}
\begin{proof}
	Differentiate \eqref{semidiscreteform} in time to see that  
\beqs
\begin{aligned}
\intd{ F(|\m_h|) \m_{ht}+ F'(|\m_h|)\frac{\m_h\cdot \m_{ht}}{|\m_h|} \m_h+ \sum_{i=-1}^N F_{a_i}(|\m_h|) a_{it} \m_h, \vv} 
-(  \brho_{ht}, \nabla\cdot \vv)= -\intd{ \nabla \Psi_t,\vv},\\
\intd{\phi \brho_{htt}, q} + \intd{ \nablad \m_{ht}, q} =\intd{ f_t, q}-\intd{\phi\Psi_{tt}, q}.
\end{aligned}
\eeqs
Taking $(\vv,q)=(\m_{ht}, \brho_{ht})$ and adding two resultant equations we obtain 
\begin{multline}\label{x1}
\frac 1 2\ddt \norm{\brho_{ht}}_{\phi}^2+\norm{F^{1/2}(|\m_h|) \m_{ht}}_{L^2}^2
 =-\intd{ F'(|\m_h|)\frac{\m_h\cdot \m_{ht}}{|\m_h|} \m_h, \m_{ht}}\\
 -\intd{ \sum_{i=-1}^N F_{a_i}(|\m_h|) a_{it} \m_h, \m_{ht}} 
 -\intd{ \nabla \Psi_t,\m_{ht}}+\intd{ f_t, \brho_{ht}} -\intd{\phi\Psi_{tt},\brho_{ht} } .
\end{multline}
We estimate the right hand side term by term 

By \eqref{dervF}, \eqref{Ka}, H\"older's inequality and Young's inequality we have 
\begin{multline}\label{x2}
-\intd{ F'(|\m_h|)\frac{\m_h\cdot \m_{ht}}{|\m_h|} \m_h, \m_{ht}}	 -\intd{ \sum_{i=-1}^N F_{a_i}(|\m_h|) a_{it} \m_h, \m_{ht}}
\le \alpha\norm{F^{1/2}(|\m_h|)\m_{ht}}_{L^2}^2 +C\intd{F(|\m_h|)|\m_h|, |\m_{ht}|}\\
\le (\alpha +\varep) \norm{F^{1/2}(|\m_h|)\m_{ht}}_{L^2}^2+ C_\varep\norm{F^{1/2}(|\m_h|)\m_h}_{L^2}^2,
\end{multline}
 and 
 \beq\label{x3}
\begin{aligned}
&-\intd{ \nabla \Psi_t,\m_{ht}}+\intd{ f_t, \brho_{ht}} -\intd{\phi\Psi_{tt},\brho_{ht} } \\
&\qquad\qquad \le \norm{f_t}_{\phi^{-1}}^2+\norm{\Psi_{tt}}_{\phi}^2+\frac 1 4\norm{\rho_{ht}}_{\phi}^2
+\varep\norm{F^{1/2}(|\m_h|)\m_{ht}}_{L^2}^2 + C_\varep\norm{F^{-1/2}(|\m_h|)\nabla\Psi_t}_{L^2}^2.
\end{aligned}
\eeq
Comparing \eqref{x1}, \eqref{x2}, \eqref{x3} and taking $\varep =(1-\alpha)/4$ yield   
\begin{align*}
	&\ddt \norm{\brho_{ht}}_{\phi}^2+(1-\alpha)\norm{F^{1/2}(|\m_h|) \m_{ht}}_{L^2}^2\\
	&\qquad\le  \norm{f_t}_{\phi^{-1}}^2+\norm{\Psi_{tt}}_{\phi}^2+\frac 1 2\norm{\brho_{ht}}_{\phi}^2
	+ C\left(\norm{F^{1/2}(|\m_h|)\m_h}_{L^2}^2  
	+ \norm{F^{-1/2}(|\m_h|)\nabla\Psi_t}_{L^2}^2\right). 
\end{align*}
By the virtual of \eqref{OrdF}, \eqref{bi2},  one has 
\beqs
\begin{aligned}
\norm{F^{-1/2}(|\m_h|)\nabla\Psi_t}_{L^2}^2+\norm{F^{1/2}(|\m_h|)\m_h}_{L^2}^2
&\le a_*^{-1}\norm{\nabla\Psi_t}_{L^2}^2+N a^*\int_\Omega (|\m_h|^{-\alpha}+|\m_h|^{\alpha_N}) \m_h^2dx\\
&\le a_*^{-1}\norm{\nabla\Psi_t}_{L^2}^2+2N a^*\int_\Omega (1+|\m_h|^{\alpha_N+2}) dx\\
&\le a_*^{-1}\norm{\nabla\Psi_t}_{L^2}^2+2N a^*(|\Omega|+1)(1+\norm{\m_h}_{L^s}^s).
\end{aligned}
\eeqs
It follows that 
\beqs
\ddt \norm{\brho_{ht}}_{\phi}^2+\norm{F^{1/2}(|\m_h|) \m_{ht}}_{L^2}^2
\le  \frac 1 2\norm{\brho_{ht}}_{\phi}^2+ C\big( \norm{f_t}_{\phi^{-1}}^2+\norm{\Psi_{tt}}_{\phi}^2+  \norm{\m_h}_{L^s}^s+\norm{\nabla\Psi_t}_{L^2}^2  \big). 
\eeqs
By Gronwall's inequality  and the boundedness of $\phi$ we find that
\begin{align*}
\norm{\brho_{ht}}_{L^\infty(I,L^2)}^2+\norm{F^{1/2}(|\m_h|) \m_{ht}}_{L^2(I,L^2)}^2 
\le C\left(\norm{\brho_{ht}(0)}_{L^2}^2 
+\norm{f_t}_{L^2(I,L^2)}^2+\norm{\Psi_{tt}}_{L^2(I,L^2)}^2+  
\norm{\m_h}_{L^s(I,L^s)}^s+\norm{\nabla\Psi_t}_{L^2(I,L^2)}^2\right).
\end{align*}
Note that by \eqref{OrdF}, one has 
\beqs
\norm{F^{1/2}(|\m_h|) \m_{ht}}_{L^2(I,L^2)}^2 \ge 
a_* \norm{\m_{ht}}_{L^2(I,L^2)}^2.
\eeqs
 Combining this fact with estimate \eqref{stab1}, we obtain the first part of estimate \eqref{stab2}. 

To verify the last part of \eqref{stab2}, we choose $q=\nablad\m_h$ in \eqref{semidiscreteform} yielding  
\beqs
\norm{\nablad\m_h}_{L^2}^2=-\intd{\phi\brho_{h,t},\nablad\m_h}+\intd{f,\nablad\m_h}-\intd{\phi\Psi_t,\nablad\m_h}.
\eeqs
Then by H\"older's inequality 
\beqs
\norm{\nablad\m_h}_{L^2}\le \norm{\phi\brho_{ht}}_{L^2}+\norm{f}_{L^2}+\norm{\phi\Psi_t}_{L^2} \le (\phi^*+1)(\norm{\brho_{ht}}_{L^2}+\norm{f}_{L^2}+\norm{\Psi_t}_{L^2}),
\eeqs
which implies
\beqs
\norm{\nablad\m_h}_{L^2}^2\le 2(\phi^*+1)^2 (\norm{\brho_{ht}}_{L^2}^2+\norm{f}_{L^2}^2+\norm{\Psi_t}_{L^2}^2).
\eeqs
Using the first part of \eqref{stab2} to bound $\norm{\brho_{ht}}_{L^2}^2$ we find that $\norm{\nablad\m_h}_{L^2}^2$ holds \eqref{stab2}. 
\end{proof}

In the same manner to the problem~\eqref{WeakProb}, we have as the following:  
 
\begin{theorem}\label{Est4sol} Suppose $(\m,\brho)$ be a solution of the problem \eqref{WeakProb}. Then, there exists  a positive constant $C$ such that
\begin{align}
\label{stab3}
&\norm{\brho}_{L^\infty(I,L^2(\Omega))}^2+\norm{\m}_{L^s(I,L^s(\Omega))}^s \le \norm{\brho_0}_{L^2(\Omega)}^2 + C\mathscr{A}, \\
\label{stab4}
&\norm{\brho_t}_{L^\infty(I,L^2(\Omega))}^2+\norm{\m_t}_{L^2(I,L^2(\Omega))}^2 +\norm{\nablad \m}_{L^\infty(I,L^2(\Omega))}^2\le C (\mathscr{B}_0 +\mathscr{B}),
\end{align} 
where $\mathscr A, \mathscr B_0, \mathscr B $ are defined in \eqref{Adef} and \eqref{Bdef}.
\end{theorem}

\section{Dependence of solutions on parameters} \label{dependence-sec}
In this section, we study the dependence of the solution on the coefficients of Forchheimer polynomial $F(s)$ in \eqref{GPoly}.
Let $N\ge 1$, the exponent vector $\alp=(\alpha_{-1},0,\alpha_1,\ldots,\alpha_N)$ and the boundary data $\psi(\x,t)$  be fixed.  Let ${\bf D}$ be a compact subset of $\{\bfa=(a_{-1},a_0,\ldots,a_N):0<a_*\le a_{-1},a_0, a_N\le a^*, 0\le a_1,\ldots,a_{N-1}\le a^*\}$.
Let $F_1(y)=F(\bfa_1,y)$ and $F_2(y)=F(\bfa_2,y)$ be two functions of class P($N,\alp$), where  $\bfa_1$ and $\bfa_2$ belong to  ${\bf D}$. Let $ \brho_h^{(1)}= \brho_h^{(1)}(\x,t;\bfa_1)$, $ \brho_h^{(2)}= \brho_h^{(2)}(\x,t;\bfa_2)$  be the two solutions of \eqref{semidiscreteform} respective to  $F(\bfa_1,y)$, $F(\bfa_2,y)$  with the same boundary data $\psi$ and initial data $\rho_0$. We will estimate $\norm{ \brho_h^{(1)}- \brho_h^{(2)}}$, $\norm{\m_h^{(1)} -\m_h^{(2)}}$ in the term of $\|\bfa_1-\bfa_2\|_{L^\infty(I,L^\infty)}$. Let $\brho_h= \brho_h^{(1)}- \brho_h^{(2)}=\rho_h^{(1)}-\rho_h^{(2)}$, $\m_h=\m_h^{(1)} - \m_h^{(2)}$. Then 
\beq\label{ss1}
\begin{split}
\intd{ F(\bfa_1,|\m_h^{(1)}|) \m_h^{(1)}-F(\bfa_2,|\m_h^{(2)}|) \m_h^{(2)}, \vv} -(  \brho_h, \nabla\cdot \vv)= 0 \quad \text{ for all } \vv\in V_h,\\
\intd{\phi \brho_{ht}, q} + \intd{ \nablad \m_h, q} =0 \quad\text{ for all } q\in Q_h. 
\end{split}
\eeq 
\begin{theorem}\label{DepCoeff}  Let $(\m_h^{(i)},\brho_h^{(i)}), i=1,2$ be two solutions to problem \eqref{semidiscreteform} corresponding to vector coefficients  $\bfa_i$ of polynominal $F(\bfa_i,y)$ in \eqref{GPoly}. There exists a positive constant $C$ independence of $\|\bfa_1-\bfa_2\|_{L^\infty(I,L^\infty(\Omega))}$  such that    
\beqs
\norm{ \brho_h^{(1)}- \brho_h^{(2)}}_{L^\infty(I,L^2(\Omega))}^2 +\norm{\m_h^{(1)}-\m_h^{(2)}}_{L^s(I,L^s(\Omega))}^s \le C \norm{\bfa_1-\bfa_2}_{L^\infty(I,L^\infty(\Omega))}^{s^*} .
\eeqs
\end{theorem}
\begin{proof} 
Choosing $(\vv,q)=(\m_h,\brho_h)$ in \eqref{ss1}, adding the resultant equations, we find that
\beq\label{Diffeq}
\frac 12\frac{d}{dt}\norm{\brho_h}_{\phi}^2+ \intd{ F(\bfa_1,|\m_h^{(1)}|) \m_h^{(1)}-F(\bfa_2,|\m_h^{(2)}|) \m_h^{(2)}, \m_h^{(1)}-\m_h^{(2)}}=0.
\eeq
According to \eqref{quasimonotone}, 
\beq\label{keyEst}
\begin{aligned}
\frac 12\frac{d}{dt}\norm{\brho_h}_{\phi}^2
\le -C_3\norm{\m_h}_{L^s}^s  dx +C_4 \int_\Omega  |\bfa_1 -\bfa_2| (1+|\m_h^{(1)}|+ |\m_h^{(2)}|)^{s-1}|\m_h|dx.
\end{aligned}
\eeq
Then by the Young's inequality 
 \begin{align*}
 	\frac{d}{dt}\norm{\brho_h}_{\phi}^2 +C_3\norm{\m_h}_{L^s}^s
 	\le C\int_\Omega  |\bfa_1 -\bfa_2|^{s^*}(1+|\m_h^{(1)}|+ |\m_h^{(2)}|)^s dx.
 \end{align*}
 Integrating in $t$ gives
 \begin{align*}
 \norm{\brho_h}_{\phi}^2 +\norm{\m_h}_{L^s(I,L^s)}^s  
 \le C\norm{\bfa_1 -\bfa_2}_{L^\infty(I,L^\infty)}^{s^*}\left(1+\norm{\m_h^{(1)}}_{L^s(I,L^s)}^s+ \norm{\m_h^{(2)}}_{L^s(I,L^s)}^s\right).
 \end{align*}
 Using \eqref{stab1}, the result \eqref{ssc2} follows. The proof is complete. 
\end{proof}

In the same manner to the problem\eqref{WeakProb}, we have as the following:
\begin{theorem} Let $(\m^{(i)},\brho^{(i)}), i=1,2$ be two solutions to problem \eqref{WeakProb} corresponding to vector coefficients  $\bfa_i$ of polynominal $F(\bfa_i,y)$ in \eqref{GPoly}. There exists a positive constant $C>0$ such that    
\beq\label{ssc2}
\norm{ \brho^{(1)}- \brho^{(2)}}_{L^\infty(I,L^2(\Omega))}^2 +\norm{\m^{(1)}-\m^{(2)}}_{L^s(I,L^s(\Omega))}^s \le C \|\bfa_1-\bfa_2\|_{L^\infty(I,L^\infty(\Omega))}^{s^*} .
\eeq
\end{theorem}
\section {Error estimates for semidiscrete approximation}\label{err-sec}
In this section, we will give the error estimate between the analytical solution and approximate solution. 
  We define the new variables:
  \beq\label{spliterr}
  \begin{split}
  \m-\m_h = \m-\Pi\m - ( \m_h- \Pi\m )=\eta - \zeta_h,\\ 
  \brho-\brho_h = \brho-\pi\brho - (\brho_h  - \pi\brho)= \theta -\vartheta_h.      
  \end{split}
  \eeq 
\begin{theorem}\label{err-thm} Let $(\m, \brho)$ be the solution of \eqref{WeakProb} and $(\m_h,\brho_h)$ be the solution of \eqref{semidiscreteform}. Suppose that $(\m, \brho)\in V \times Q$, and  $\brho_t\in L^2(I,L^2(\Omega) )$. Then there exists a positive constant $C$ independence of $h$ such that    
\beq\label{eest1}
\begin{aligned}
	\norm{\brho-\brho_h}_{L^\infty(I,L^2(\Omega))}^2+\norm{\m-\m_h}_{L^s(I,L^s(\Omega))}^s \le C h^{s^*(1-\alpha)}. 
	\end{aligned}
\eeq
\end{theorem}
\begin{proof}
Error equations 
\beq \label{ErrEq}
\begin{aligned}
	\intd{ F(|\m |) \m - F(|\m_h|) \m_h, \vv} -\intd{ \brho- \brho_h, \nabla\cdot \vv}= 0 \quad \text{ for all } \vv\in V_h,\\
	\intd{\phi (\brho_{t}-\brho_{ht}), q} + \intd{\nablad (\m-\m_h) , q} =0 \quad\text{ for all } q\in Q_h.
\end{aligned}
\eeq
Using $L^2$-project and Raviar-Thomas projection we rewrite equation as form  
\begin{align*}
	\intd{ F(|\m |) \m - F(|\Pi\m|) \Pi\m, \vv} +\intd{ F(|\Pi\m|) \Pi\m - F(|\m_h|) \m_h, \vv} 
	+ \intd{ \vartheta_h, \nabla\cdot \vv}= 0,\\
	\intd{\phi \theta_t, q} -\intd{\phi \vartheta_{ht}, q} -\intd{\nablad \zeta_h , q} =0.
\end{align*}
Take $q=-\vartheta_h\in Q_h$ and $\vv=-\zeta_h\in V_h  $. Add  these two equations together 
\beqs
\intd{\phi \vartheta_{ht}, \vartheta_h}+\intd{ F(|\Pi\m|) \Pi\m - F(|\m_h|) \m_h, \Pi\m-\m_h}
=	\intd{ F(|\m |) \m - F(|\Pi\m|) \Pi\m, \zeta_h}
    +\intd{\phi \theta_t, \vartheta_h}.
\eeqs
Then by the \eqref{monotone0},  
\beqs
\frac 1 2 \ddt\norm{\vartheta_h}_{\phi}^2+C_3\norm{\zeta_h}_{L^s}^s
\le \intd{ F(|\m |) \m - F(|\Pi\m|) \Pi\m, \zeta_h}+\intd{\phi \theta_t, \vartheta_h}.
\eeqs
Using Young's and H\"older's inequality we find that  
\beqs
\intd{\phi \theta_t, \vartheta_h} \le (4\varep)^{-1}\norm{\theta_t}_{\phi}^2+\varep \norm{\vartheta_h}_{\phi}^2,
\eeqs
and 
\beq\label{sas}
\begin{split}
 \intd{ F(|\m |) \m - F(|\Pi\m|) \Pi\m, \zeta_h} 
 &\le \intd{|F(|\m |) \m - F(|\Pi\m|) \Pi\m|,|\zeta_h|}\\
 & \le  C_1\intd {( 1+|\Pi\m|+|\m|)^{s-2+\alpha}|\eta|^{1-\alpha}, |\zeta_h|} \quad \text{ (by \eqref{Lipchitz}) }\\
 &\le \varep \norm{\zeta_h}_{L^s}^s +C_\varep \norm{ ( 1+|\Pi\m|+|\m|)^{s-2+\alpha}|\eta|^{1-\alpha} }_{L^{s^*}}^{s^*} \\
 &\le \varep \norm{\zeta_h}_{L^s}^s +C_\varep \left( 1+\norm{\Pi\m}_{L^{s}}^{s}+\norm{\m}_{L^{s}}^{s}\right)^{\frac{s-2+\alpha}{s-1} } \norm{\eta}_{L^s}^{s^*(1-\alpha)}.
\end{split}
\eeq
Choose $\varep=C_3/2$, we find that 
\beqs
\ddt\norm{\vartheta_h}_{\phi}^2+C_3\norm{\zeta_h}_{L^s}^s \le \varep \norm{\vartheta_h}_{\phi}^2
+C\left( 1+\norm{\Pi\m}_{L^{s}}^{s}+\norm{\m}_{L^{s}}^{s}\right)^{\frac{s-2+\alpha}{s-1} }\norm{\eta}_{L^s}^{s^*(1-\alpha)}+(4\varep)^{-1}\norm{\theta_t}_{\phi}^2.
\eeqs

Integrating in time from $0$ to $t$ and take sup-norm
\begin{align*}
\sup_{t\in[0,T]}\norm{\vartheta_h}_{\phi}^2+\norm{\zeta_h}_{L^s(I,L^s)}^s
&\le \varep T \sup_{t\in[0,T]} \norm{\vartheta_h}^2+(4\varep)^{-1}\int_0^T \norm{\theta_t}_{\phi}^2 dt\\
&\quad+C\int_0^T \left( 1+\norm{\Pi\m}_{L^{s}}^{s}+\norm{\m}_{L^{s}}^{s}\right)^{\frac{s-2+\alpha}{s-1} } \norm{\eta}_{L^s}^{s^*(1-\alpha)}dt.
\end{align*}
Now take $\varep =1/(2T)$, we find that
\beqs
\sup_{t\in[0,T]}\norm{\vartheta_h}_{\phi}^2+\norm{\zeta_h}_{L^s(I,L^s)}^s
\le \frac T 2 \int_0^T \norm{\theta_t}_{\phi}^2 dt
+C\int_0^T \left( 1+\norm{\Pi\m}_{L^{s}}^{s}+\norm{\m}_{L^{s}}^{s}\right)^{\frac{s-2+\alpha}{s-1} } \norm{\eta}_{L^s}^{s^*(1-\alpha)}dt.
\eeqs
Thus 
\beq\label{d1}
\norm{\vartheta_h}_{L^\infty(I,L^2)}^2+\norm{\zeta_h}_{L^s(I,L^s)}^s
\le \frac T 2 \int_0^T \norm{\theta_t}_{\phi}^2 dt
+C\int_0^T \left( 1+\norm{\Pi\m}_{L^{s}}^{s}+\norm{\m}_{L^{s}}^{s}\right)^{\frac{s-2+\alpha}{s-1} } \norm{\eta}_{L^s}^{s^*(1-\alpha)}dt.
\eeq
Using \eqref{bi0}, \eqref{d1} we obtain 
 \begin{align*}
  &\norm{\brho-\brho_h}_{L^\infty(I,L^2)}^2+\norm{\m-\m_h}_{L^s(I,L^s)}^s\\
  &\quad \le 2\left(\norm{\theta}_{L^\infty(I,L^2)}^2+ \norm{\vartheta_h}_{L^\infty(I,L^2)}^2\right) +2^{s-1}\left(\norm{\eta}_{L^s(I,L^s)}^s+ \norm{\zeta_h}_{L^s(I,L^s)}^s\right)\\
   &\quad \le 2^{s-1}\left(\norm{\theta}_{L^\infty(I,L^2)}^2 +\norm{\eta}_{L^s(I,L^s)}^s+\norm{\vartheta_h}_{L^\infty(I,L^2)}^2+\norm{\zeta_h}_{L^s(I,L^s)}^s  \right)\\
   &\quad \le 2^{s-1}\Big(\norm{\theta}_{L^\infty(I,L^2)}^2 +\norm{\eta}_{L^s(I,L^s)}^s 
  + T  \int_0^T \norm{\theta_t}_{\phi}^2 dt
+C\int_0^T \left( 1+\norm{\Pi\m}_{L^{s}}^{s}+\norm{\m}_{L^{s}}^{s}\right)^{\frac{s-2+\alpha}{s-1} } \norm{\eta}_{L^s}^{s^*(1-\alpha)}dt\Big).
\end{align*}
Applying estimates \eqref{divProj} and \eqref{L2Proj} implies that 
\begin{multline*}
\norm{\brho-\brho_h}_{L^\infty(I,L^2)}^2+\norm{\m-\m_h}_{L^s(I,L^s)}^s
\le Ch^2\Big(\norm{\brho}_{L^\infty(I,L^2)}^2+ \norm{\brho_t}_{L^2(I,L^2)}^2\Big)\\
+C h^s \norm{\m}_{L^s(I,W^{1,s})}^s+   Ch^{s^*(1-\alpha)}\int_0^T \left(1+\norm{\Pi\m}_{L^{s}}^{\frac{s(s-2+\alpha)}{s-1}  }+\norm{\m}_{L^{s}}^{\frac{s(s-2+\alpha)}{s-1}} \right) \norm{\m}_{W^{1,s}}^{s^*(1-\alpha)} dt.
\end{multline*}

Note that $\frac{s(s-2+\alpha)}{s-1} + s^*(1-\alpha) =s$,
\begin{align*}
&\left(1+\norm{\Pi\m}_{L^{s}}^{\frac{s(s-2+\alpha)}{s-1}  }+\norm{\m}_{L^{s}}^{\frac{s(s-2+\alpha)}{s-1}} \right) \norm{\m}_{W^{1,s}}^{s^*(1-\alpha)}\\
&\quad \le  C\left(1+\norm{\m}_{W^{1,s}}^{s}\right)
               +C \left(\norm{\Pi\m}_{L^{s}}^{s}+\norm{\m}_{W^{1,s}}^{s} \right)
               +C \left(\norm{\m}_{L^{s}}^{s}+\norm{\m}_{W^{1,s}}^{s} \right) \quad \text {(Young's inequality)}\\
&\quad \le 
            C \left(1+\norm{\Pi\m}_{L^{s}}^{s}+\norm{\m}_{W^{1,s}}^{s} \right)\quad \text{(By  \eqref{DProj})}\\
            &\quad \le 
            C \left(1+\norm{\m}_{L^{s}}^{s}+h\norm{\nablad \m}_{L^2}^{s}+\norm{\m}_{W^{1,s}}^{s} \right).
\end{align*}

Thus 
\begin{align*}
\norm{\brho-\brho_h}_{L^\infty(I,L^2)}^2+\norm{\m-\m_h}_{L^s(I,L^s)}^s
&\le Ch^2\Big(\norm{\brho}_{L^\infty(I,L^2)}^2+ \norm{\brho_t}_{L^2(I,L^2)}^2\Big)\\
&+C h^s \norm{\m}_{L^s(I,W^{1,s})}^s+   Ch^{s^*(1-\alpha)}\int_0^T (1+\norm{\Pi\m}_{L^{s}}^s+\norm{\m}_{W^{1,s}}^{s})  dt.
\end{align*}
Since $s^*(1-\alpha)<2<s,$ 
\beqs
\norm{\brho-\brho_h}_{L^\infty(I,L^2)}^2+\norm{\m-\m_h}_{L^s(I,L^s)}^s
 \le Ch^{s^*(1-\alpha)} \Big( 1+ \norm{\brho}_{L^\infty(I,L^2)}^2+ \norm{\brho_t}_{L^2(I,L^2)}^2
 +  \norm{\m}_{L^s(I,W^{1,s})}^s+\norm{\Pi\m}_{L^s(I,L^s)}^s  \Big).
\eeqs
The proof is concluded.  
\end{proof}
\begin{theorem}\label{par-thm}
Let $(\m_h^{(i)}, \brho_h^{(i)}), i=1,2$ be two solutions to problems \eqref{semidiscreteform} corresponding to vector coefficients  $\bfa_i$ of generalized polynominals $F(\bfa_i,s)$ in \eqref{GPoly}. Suppose that $(\m^{(i)}, \brho^{(i)})\in V\times Q$  and $\rho_t^{(i)}\in L^2(I,L^2(\Omega) )$. Then, there exists a constant positive constant $C$ independent of $h$ and $\norm{\bfa_1 -\bfa_2}_{L^\infty(I,L^\infty)}$ such that   
\beq\label{mrhoDif}
\norm{ \brho_h^{(1)} - \brho_h^{(2)}}_{L^\infty(I,L^2(\Omega))}^2+ \norm{\m_h^{(1)} - \m_h^{(2)}}_{L^s(I,L^s(\Omega))}^s \le C \big( h^{s^*(1-\alpha)} + \norm{\bfa_1 -\bfa_2}_{L^\infty(I,L^\infty(\Omega))}^{s^* } \big).
\eeq
\end{theorem}
\begin{proof}
The triangle inequality shows that 
\begin{align*}
\norm{ \brho_h^{(1)} - \brho_h^{(2)}}_{L^\infty(I,L^2)}^2 + \norm{\m_h^{(1)} -\m_h^{(2)}}_{L^s(I,L^s)}^s 
&\le C \sum_{i=1}^2\Big(\norm{\brho_h^{(i)} -\brho^{(i)}}_{L^\infty(I,L^2)} + \norm{\m_h^{(i)} -\m^{(i)}}_{L^\infty(I,L^\infty)}^2\Big)\\
& +C\Big(\norm{\rho^{(1)} -\rho^{(2)} }_{L^\infty(I,L^2)}^2+\norm{\m^{(1)} -\m^{(2)}}_{L^s(I,L^s)}^s\Big).
\end{align*}
Then by using \eqref{eest1} to treat the sum-term and \eqref{ssc2} to the last terms we obtain \eqref{mrhoDif}.   
\end{proof}
\section{Fully discrete method}\label{err-ful-sec}
  Let $\{t_n\}_{n=0}^M$ be the uniform partition of $[0,T]$ with $t_n=n\Dt$, for time step $\Dt >0$. We define $\varphi_n =\varphi(\cdot, t_n) $.  The discrete time mixed finite element approximation to \eqref{semidiscreteform} is defined as follows:  

For given $\brho_{h0}(x)=\pi \brho_0(x)$ and $\big\{f_n\big\}_{n=1}^M\in L^2(\Omega), \big\{\Psi_n\big\}_{n=1}^M, \big\{\Psi_{tn}\big\}_{n=1}^M\in C(\bar\Omega)  $. 
Find a pair $ (\m_{hn}, \brho_{hn} )$ in  $V_h \times Q_h $, $n=0,1, 2,\ldots, M$  such that 

\beq\label{fullydiscreteform} 
\begin{aligned}
  \intd{ F(|\m_{hn}|)\m_{hn}, \vv} - \intd{\brho_{hn}, \nabla\cdot \vv} = -\intd{ \nabla\Psi_n,\vv } \quad\text{ for all } \vv\in V_h,\\
\intd{\phi \frac {\brho_{hn}-\brho_{hn-1}}{\Dt}, q} + \intd{\nablad \m_{hn} , q} =\intd{ f_n , q} -\intd{\phi\Psi_{tn}, q}  \quad\text{ for all } q\in Q_h,
\end{aligned}
\eeq 
for $\brho_{h0}(x)=\pi \brho_0(x)$. 
\begin{lemma}[Stability]\label{stab-appr} Let  $( \m_{hn}, \brho_{hn})$ solve the fully discrete finite element
	approximation \eqref{fullydiscreteform} for each time step $n=1,2,\ldots ,M$. There exists a positive constant $C$ independent of $n,\Dt$ such that for $\Dt$ sufficiently small 
	\beq\label{pwBound}
	\norm{\brho_{hn} }_{L^2(\Omega)}^2 +\sum_{i=1}^n \Dt\norm{\m_{hi}}_{L^s(\Omega)}^s \le C   \Big(\norm{\brho_0}_{L^2(\Omega)}^2 + \sum_{i=1}^n  \Dt (\norm{f_n}_{L^2(\Omega)}^2 + \norm{\Psi_{tn}}_{L^2(\Omega)}^2 + \norm{\nabla\Psi_n}_{L^{s^*}(\Omega)}^{s^*})\Big) .  
	\eeq
\end{lemma}
\begin{proof}
	Selecting $(\vv,q)=2(\m_{hn}, \brho_{hn})$ in \eqref{fullydiscreteform}, we find that  
	\beq
	\begin{split}
		 2\intd{ F(|\m_{hn}|)\m_{hn}, \m_{hn}} - 2\intd{\brho_{hn}, \nabla\cdot \m_{hn}} = -2\intd{\nabla\Psi_n ,\m_{hn}},\\
		2\intd{\phi \frac {\brho_{hn}-\brho_{hn-1}}{\Dt}, \brho_{hn}} + 2\intd{\nablad \m_{hn} , \brho_{hn}} =2\intd{ f_n ,\brho_{hn}}-2\intd{\phi\Psi_{tn},\brho_{hn}}. 
	\end{split}
	\eeq
	Adding the two above equations, and using the identity  
	\beqs
	2\Big( \phi (\brho_{hn} - \brho_{hn-1}), \brho_{hn}\Big) = \norm{\brho_{hn}}_{\phi}^2 - \norm{\brho_{hn-1}}_{\phi}^2 +\norm{\brho_{hn} -\brho_{hn-1}}_{\phi}^2,
	\eeqs 
	we obtain 
	\begin{multline}\label{t1}
	\norm{\brho_{hn}}_{\phi}^2 - \norm{\brho_{hn-1}}_{\phi}^2 +\norm{\brho_{hn} -\brho_{hn-1}}_{\phi}^2 + 2\Dt\intd{ F(|\m_{hn}|)\m_{hn}, \m_{hn}}
	= 2\Dt \Big(\intd{ f_n ,\brho_{hn}}-\intd{\phi\Psi_{tn},\brho_{hn}}-\intd{\nabla\Psi_n ,\m_{hn}}\Big).  
	\end{multline}
	It follows from \eqref{monotone0} that
	\beq\label{t2}
	\intd{ F(|\m_{hn}|)\m_{hn}, \m_{hn}} \ge C_3\norm{\m_{hn}}_{L^s}^s. 
	\eeq
	Using H\"older's inequality to the RHS of \eqref{t1} shows that
	\begin{multline}\label{t3}
	 \intd{ f_n ,\brho_{hn}}-\intd{\phi\Psi_{tn},\brho_{hn}}-\intd{\nabla\Psi_n ,\m_{hn}}
	\le \norm{f_n}_{\phi^{-1}}^2 + \norm{\Psi_{tn}}_{\phi}^2 + \frac {1}{2}\norm{\brho_{hn}}_{\phi}^2+\frac{C_3}{2} \norm{\m_{hn}}_{L^s}^s + \frac{(sC_3/2)^{-s^*/s}}{s^*}  \norm{\nabla\Psi_n}_{L^{s^*}}^{s^*}.
	\end{multline}
	Combining \eqref{t1}--\eqref{t3}  yields  
	\beqs
	\begin{split}
		\norm{\brho_{hn}}_{\phi}^2 - \norm{\brho_{hn-1}}_{\phi}^2  +C_3\Dt \norm{\m_{hn}}_{L^s}^{s}
		&\le  \Dt\norm{\brho_{hn}}_{\phi}^2 +C\Dt\left(\norm{f_n}_{\phi^{-1}}^2 + \norm{\Psi_{tn}}_{\phi}^2 + \norm{\nabla\Psi_n}_{L^{s^*}}^{s^*}\right),
	\end{split}
	\eeqs
	which shows that    
	\beqs
	\frac{\norm{\brho_{hn}}_{\phi}^2 - \norm{\brho_{hn-1}}_{\phi}^2}{\Dt} - \norm{\brho_{hn}}_{\phi}^2 +C_3 \norm{\m_{hn}}_{L^s}^{s} \le   C\left(\norm{f_n}_{\phi^{-1}}^2 + \norm{\Psi_{tn}}_{\phi}^2 + \norm{\nabla\Psi_n}_{L^{s^*}}^{s^*}\right).
	\eeqs
		
	By discrete Gronwall's inequality in Lemma~\ref{DGronwall},  
	\beq\label{dscEst1}
	\begin{split}
	\norm{\brho_{hn} }_{\phi}^2 + C_3\sum_{i=1}^n \Dt  \norm{\m_{hi}}_{L^s}^{s}  \le 
	 C e^{\frac{n\Dt}{1-\Dt}}  \norm{\brho_{h0}}_{\phi}^2 +Ce^{\frac{n\Dt}{1-\Dt}}\sum_{i=1}^n \Dt\left(\norm{f_i}_{\phi^{-1}}^2 + \norm{\Psi_{ti}}_{\phi}^2 + \norm{\nabla\Psi_i}_{L^{s^*}}^{s^*}\right).
	 \end{split}
	\eeq
Note that $\norm{\brho_{h0}}_{L^2}^2 \le \norm{\brho_0}_{L^2}^2$  and $ e^{\frac{n\Dt}{1-\Dt}} \le e^{\frac{N\Dt}{1-\Dt}} = e^{\frac{NT}{N-T}} $ imply \eqref{pwBound}. The proof is complete.
\end{proof}
\subsection{Error analysis}
 
As in the semidiscrete case, we use $\eta = \m -\Pi \m,$ $\zeta_h=\m_h-\Pi \m $,   $\theta=\brho-\pi\brho$, $\vartheta_h=\brho_h-\pi \brho$ and $\eta_n, \theta_n$, $\zeta_{hn},\vartheta_{hn}$  be evaluating $\eta, \theta$, $\zeta_h,\vartheta_h$  at the discrete time levels.
We also define  
$$\partial \varphi_n = \frac {\varphi_n -\varphi_{n-1} }{\Dt}.$$
\begin{theorem}\label{Err-ful}
	Let $(\m_n,\brho_n)$ solve problem \eqref{semidiscreteform} and $(\m_{hn},\brho_{hn})$ solve the fully  discrete finite element approximation \eqref{fullydiscreteform} for each time step $n$, $n=1,\ldots, M$. Suppose that  $(\m,\brho)\in V\times Q  $ and $\brho_{tt} \in L^2(I,L^2(\Omega))$. Then, there exists a positive constant $C$ independent of $h$ and $\Dt$  such that for $\Dt$ sufficiently small  
	\beq\label{fulerrl2}
	\norm {\brho_n-\brho_{hn}}_{L^2(\Omega)}^2 +\norm {\m_n -\m_{hn}}_{L^s(\Omega)}^s  \le C\big( h^{s^*(1-\alpha)}+ \Dt^2 \big).
	\eeq  
\end{theorem}
\begin{proof}
 Evaluating equation \eqref{WeakProb} at $t=t_n$ gives
  \beq\label{fuldis1}
 \begin{aligned}
 	\intd{ F(|\m_n|) \m_n, \vv} -(  \brho_n, \nabla\cdot \vv)= -\intd{\nabla\Psi_n ,\vv} \quad \text{ for all } \vv\in V_h,\\
 	\intd{\phi \brho_{tn}, q} + \intd{ \nablad \m_n, q} =\intd{ f_n ,q}-\intd{\phi\Psi_{tn},q} \quad\text{ for all } q\in Q_h.
 \end{aligned}
 \eeq

 Subtracting \eqref{fullydiscreteform} from \eqref{fuldis1}, we obtain 
 \begin{align}
  \intd{ F(|\m_n|) \m_n- F(|\m_{hn}|)\m_{hn} , \vv} -(  \pi\brho_n-\brho_{hn}, \nabla\cdot \vv)= 0, \quad\text{ for all } \vv\in V_h,\\
    \intd{\phi(\brho_{tn} -  \frac {\brho_{hn}-\brho_{hn-1}}{\Dt}) , q} + \intd{\nablad (\Pi\m_n-\m_{hn}) , q} =0   \quad\text{ for all } q\in Q_h
 \end{align}
Choosing $\vv=-\zeta_{hn}, q=-\vartheta_{hn}$, and adding the two equations shows that   
\beq\label{errEq}
\intd{\phi (\brho_{tn}-\partial\brho_{hn}),\vartheta_{hn}} - \intd{ F(|\m_n|) \m_n- F(|\m_{hn}|)\m_{hn} , \Pi \m_n-\m_{hn}}  =0 . 
\eeq 
Since  $\brho_{tn}-\partial\brho_{hn} = \brho_{tn}-\partial \brho_n+ \partial \theta_n - \partial \vartheta_{hn},$  
we rewrite \eqref{errEq} in the form  
\beq\label{eereq}
\begin{split}
	&(\phi\partial \vartheta_{hn} , \vartheta_{hn} ) + \intd{ F(|\Pi\m_n|) \Pi\m_n- F(|\m_{hn}|)\m_{hn} , \Pi \m_n-\m_{hn}}\\
	&\quad= \intd{  F(|\m_n|)\m_n-F(|\Pi\m_n|) \Pi\m_n , \zeta_{hn}}+\intd{\phi (\brho_{tn}-\partial \brho_n),\vartheta_{hn} }+\intd{\phi\partial \theta_n,\vartheta_{hn}}.  
\end{split} 
\eeq
We will evaluate \eqref{eereq}  term by term. 

\textbullet\quad  For the first term, we use the identity 
\beq\label{rhs1}
	(\phi\partial \vartheta_{hn} , \vartheta_{hn} )= \frac 1{2\Dt}\left(\norm{\vartheta_{hn}}_{\phi}^2  -\norm{\vartheta_{hn-1}}_{\phi}^2\right)+\frac{\Dt}{2}\norm{\partial \vartheta_{hn} }_{\phi}^2 . 
\eeq

\textbullet\quad For the second term, the monotonicity of $F(\cdot)$ in \eqref{monotone0} yields 
\beq\label{rhs2} 
	\intd{ F(|\Pi\m_n|) \Pi\m_n- F(|\m_{hn}|)\m_{hn} , \Pi \m_n-\m_{hn}} \ge C_3 \norm{\zeta_{hn}}_{L^s}^s. 
\eeq

\textbullet\quad For the third term, using \eqref{sas} gives       
\beq\label{lhs1}
\intd{  F(|\m_n|)\m_n-F(|\Pi\m_n|) \Pi\m_n , \zeta_{hn}}
\le \frac {C_3}{2} \norm{\zeta_{hn}}_{L^s}^s +C \left( 1+\norm{\Pi\m_n}_{L^{s}}^{s}+\norm{\m_n}_{L^{s}}^{s}\right)^{\frac{s-2+\alpha}{s-1} } \norm{\eta_n}_{L^s}^{s^*(1-\alpha)}.
\eeq

\textbullet\quad Using Cauchy-Schwarz's inequality, Poincare's and Young's inequalities, we obtain    
\beq\label{lhs2}
\intd{\phi (\brho_{tn}-\partial \brho_n),\vartheta_{hn} }+\intd{\phi\partial \theta_n,\vartheta_{hn}}
	\le C\Dt \int_{t_{n-1}}^{t_n} \norm{\brho_{tt}(\tau)}_{\phi}^2 d\tau  
	+\frac{C}{\Dt} \int_{t_{n-1}}^{t_n}\norm{\theta_t(\tau)}_{\phi}^2 d\tau   +\frac {1} 2 \norm{\vartheta_{hn}}_{\phi}^2.
\eeq

In view of \eqref{rhs1}--\eqref{lhs2}, \eqref{eereq} yields       
\beq
\begin{split}
\frac {\norm{\vartheta_{hn}}_{\phi}^2  -\norm{\vartheta_{hn-1}}_{\phi}^2}{\Dt} -  \norm{\vartheta_{hn}}_{\phi}^2  + C_3\norm{\zeta_{hn}}_{L^s}^s
\le C\Dt \int_{t_{n-1}}^{t_n} \norm{\brho_{tt}(\tau)}_{\phi}^2 d\tau + \frac{C}{\Dt} \int_{t_{n-1}}^{t_n}\norm{\theta_t(\tau)}_{\phi}^2 d\tau\\
+C\left( 1+\norm{\Pi\m_n}_{L^{s}}^{s}+\norm{\m_n}_{L^{s}}^{s}\right)^{\frac{s-2+\alpha}{s-1} } \norm{\eta_n}_{L^s}^{s^*(1-\alpha)}.	
\end{split}
\eeq
By mean of discrete Gronwall's inequality in Lemma~\ref{DGronwall} and the fact $\vartheta_{h0}=0$, we find that  
\beqs
\begin{split}
	\norm{\vartheta_{hn}}_{\phi}^2 + C_3\sum_{i=1}^n \Dt\norm{\zeta_{hi}}_{L^s}^s
	&\le Ce^{\frac{n\Dt}{1-\Dt}}\Dt\sum_{i=1}^n \Dt \int_{t_{i-1}}^{t_i} \norm{\brho_{tt}(\tau)}_{\phi}^2 d\tau+ \frac{1}{\Dt} \int_{t_{i-1}}^{t_i} \norm{\theta_t(\tau)}_{\phi}^2 d\tau\\
	&+Ce^{\frac{n\Dt}{1-\Dt}}\Dt\sum_{i=1}^n \left( 1+\norm{\Pi\m_i}_{L^{s}}^{s}+\norm{\m_i}_{L^{s}}^{s}\right)^{\frac{s-2+\alpha}{s-1} } \norm{\eta_i}_{L^s}^{s^*(1-\alpha)}\\
	&\le Ce^{\frac{N\Dt}{1-\Dt}} \Dt^2 \int_0^T \norm{\brho_{tt}(\tau)}_{\phi}^2 d\tau + \int_0^T \norm{\theta_t(\tau)}_{\phi}^2 d\tau\\
	&+Ce^{\frac{N\Dt}{1-\Dt}}\Dt\sum_{i=1}^n\left( 1+\norm{\Pi\m_i}_{L^{s}}^{s}+\norm{\m_i}_{L^{s}}^{s}\right)^{\frac{s-2+\alpha}{s-1} } \norm{\eta_i}_{L^s}^{s^*(1-\alpha)}.
\end{split}
\eeqs
Then it follows 
\beq\label{comb3}
\begin{split}
	\norm{\vartheta_{hn}}_{L^2}^2 + C_3\sum_{i=1}^n \Dt\norm{\zeta_{hi}}_{L^s}^s 
	&\le C \Dt^2 \int_0^T \norm{\brho_{tt}(\tau)}_{L^2}^2 d\tau + \int_0^T \norm{\theta_t(\tau)}^2 d\tau\\
	&+C\Dt\sum_{i=1}^N\left( 1+\norm{\Pi\m_i}_{L^{s}}^{s}+\norm{\m_i}_{L^{s}}^{s}\right)^{\frac{s-2+\alpha}{s-1} } \norm{\eta_i}_{L^s}^{s^*(1-\alpha)}.
\end{split}
\eeq

Using \eqref{bi0}  and \eqref{comb3} we find that 
\begin{align*}
\norm{\brho_n-\brho_{hn}}_{L^2}^2 + \sum_{i=1}^n \Dt \norm{\m_{i} -\m_{hi} }_{L^s}^s 
 &\le 2\left(\norm{\theta_n}_{L^2}^2+ \norm{\vartheta_{hn}}_{L^2}^2\right) +2^{s-1}\sum_{i=1}^n \Dt\left(\norm{\eta_i}_{L^s}^s+ \norm{\zeta_{hi}}_{L^s}^s\right)\\
 &\le 2^{s-1}\left( \norm{\theta_n}_{L^2}^2 +\sum_{i=1}^n \Dt\norm{\eta_i}_{L^s}^s + \norm{\vartheta_{hn}}_{L^2}^2+\sum_{i=1}^n \Dt\norm{\zeta_{hi}}_{L^s}^s  \right)\\
&\le  C\left(\norm{\theta_n}_{L^2}^2 +\sum_{i=1}^n \Dt\norm{\eta_i}_{L^s}^s +  \Dt^2 \int_0^T \norm{\brho_{tt}(\tau)}_{L^2}^2 d\tau + \int_0^T \norm{\theta_t(\tau)}_{L^2}^2 d\tau\right)\\
&\quad +C\sum_{i=1}^n \Dt\left( 1+\norm{\Pi\m_i}_{L^{s}}^{s}+\norm{\m_i}_{L^{s}}^{s}\right)^{\frac{s-2+\alpha}{s-1} } \norm{\eta_i}_{L^s}^{s^*(1-\alpha)}.
\end{align*}
Applying \eqref{divProj} and \eqref{L2Proj}, \eqref{bi0} and Young's inequality gives
\beqs
\begin{aligned}
\norm{\brho_n-\brho_{hn}}_{L^2}^2 +\sum_{i=1}^n \Dt \norm{\m_{i} -\m_{hi} }_{L^s}^s 
&\le  C(h^2\norm{\brho_n}_{L^2}^2 +h^s \sum_{i=1}^n \Dt\norm{\m_i}_{L^s}^s 
+  \Dt^2 \norm{\brho_{tt}(\tau)}_{L^2(I,L^2)}^2 + h^2\norm{\brho_t}_{L^2(I,L^2)}^2)\\
&\quad+Ch^{s^*(1-\alpha)}\sum_{i=1}^n  \Dt\left(1+\norm{\Pi\m_i}_{L^{s}}^{s}+\norm{\m_i}_{W^{1,s}}^{s}\right).
\end{aligned}
\eeqs 
The proof is complete.  
\end{proof}

The following theorem about an error estimate for $(\m_{hn}, \brho_{hn})$ is obtained by using the same manner as in the proof of Theorem~\ref{DepCoeff}.     

\begin{theorem}\label{Dep-ful} Let $(\m_n^{(i)},\brho_n^{(i)})$, i=1,2 solve problem \eqref{semidiscreteform} and $(\m_{hn}^{(i)}, \brho_{hn}^{(i)})$ solve the fully  discrete finite element approximation \eqref{fullydiscreteform} corresponding to vector coefficients  $\bfa_i$ of generalized polynomials $F( \bfa_i, s)$ in \eqref{GPoly}  for each time step $n$, $n=1,\ldots, M$. Suppose that each $(\m^{(i)},\brho^{(i)})\in  V\times Q $  and $\brho_{tt}^{(i)}\in L^2(I,L^2(\Omega) )$. There exists a positive constant  $C$  independent of $h,$ $\Dt$ and $ \norm{ \bfa_1 -\bfa_2 }_{L^\infty(I,L^\infty)}$  such that 
\beq\label{ErrContEst}
\norm {\brho_{hn}^{(1)} -\brho_{hn}^{(2)}}_{L^2(\Omega)}^2 + \norm {\m_{hn}^{(1)} -\m_{hn}^{(2)}}_{L^s(\Omega)}^s \le C\left( h^{s^*(1-\alpha)}+\Dt^2 + \norm{\bfa_1 -\bfa_2}_{L^\infty(I,L^\infty(\Omega))}^{s^*} \right).
\eeq  
\end{theorem}
\section{Numerical results} \label{Num-result}
In this section we carry out numerical simulations using mixed finite element approximation to solve problem \eqref{fullydiscreteform} in two dimensions to validate our theoretical estimates. For simplicity, the region of examples are unit square $\Omega=[0,1]\times[0,1]$.  
We use the piecewise-linear elements for the both density and momentum variables. We divided the unit square into an $ \mathcal N \times \mathcal N$ mesh of squares, each of them subdivided into two right triangles. For each mesh, we solved the solve problem \eqref{fullydiscreteform} numerically. Our problem is solved at each time level starting at $t=0$ until the given final time $T=1$, $I=[0,1]$.  A Newton iteration was used to solve the nonlinear equation generated at each time step. The error control in each nonlinear solve is $\rm{tol}=10^{-6}.$ At time $T$, we measured the error in $L^2$-norm for density and $L^3$-norm for the vector momentum.  We obtain the convergence rates $r=\frac{\ln(e_i/e_{i-1})}{\ln(h_i/h_{i-1})}$ of finite approximation at 7 levels with the discretization parameters $h=1/4, 1/8,1/16, 1/32, 1/64, 1/128, 1/256$ respectively and time step $\Delta t =h/2.$ 
  The numerical examples in this section are constructed in two categories:    
   \begin{itemize}
   \item Examples 1 is used to study the convergence rates of the method proposed in the paper. 
   \item Example 2 is used to study the dependence of solution on physical parameters. 
    \end{itemize}
   
 {\bf Example 1.} We test the convergence of our method with the function $F(y)=y^{-1/2}+1+y$. In this case, $\alpha=1/2$, $\alpha_N=1$, $s=3$, $s^*=3/2$, $s^*(1-\alpha)=3/4$. To test the convergence rates, we choose the analytical solution  
\beqs
\begin{split}
\rho(\x,t)=\frac 1{\sqrt 2} \Big(e^{-t} +e^{-2t}+e^{-4t}\Big)(x_1+ x_2) , \,  \text{ and }\,
\m(\x,t) =-e^{-2t}\left(\frac 1{\sqrt 2}, \frac 1{\sqrt 2}\right) \quad \forall  (\x,t)\in \Omega\times I.
\end{split}
\eeqs

For simplicity, we take $\phi(\x)\equiv 1$ on $\Omega$. The forcing term $f$ is determined from equation $\rho_t +\nabla \cdot \m = f$. Explicitly, 
$$
f(\x,t)=-\frac 1 {\sqrt 2}\Big(e^{-t} +2e^{-2t}+ 4 e^{-4t}\Big)(x_1+ x_2).
$$

The initial condition and boundary condition are determined according to the analytical solution as follows: 
\begin{align*}
 \rho_0(x) =\frac 3{\sqrt 2}(x_1+x_2), \quad 
 \psi(\x,t)=\frac 1{\sqrt 2}\Big( e^{-t} +e^{-2t}+ e^{-4t}\Big)\begin{cases}  
x_2  & \text { on } \{0\}\times[0,1] ,\\
1+x_2  & \text { on } \{1\}\times[0,1], \\
x_1  & \text { on } [0,1]\times \{0\}, \\
1+x_1  & \text { on } [0,1] \times \{1\}.
 \end{cases}
\end{align*}

Then $\Psi(\x,t)=\rho(\x,t)$. Thus $\brho(\x,t)=0$.

We report the approximate errors  in Table 1. 

 \vspace{0.2cm}  
\begin{center}
\begin{tabular}{l| c| c| c| c}
\hline
 $\mathcal N$  & \qquad $\norm{\brho-\brho_h}_{L^2(\Omega)}$ \qquad  &  \quad Rates  \quad&  \quad$\norm{\m-\m_h}_{L^3(\Omega)}$ \quad & \quad  Rates \quad  \\
\hline
4	&$2.566e-1$        	& --         		&$4.823e-1$ 		&--\\
8	&$1.689e-1$        	& $0.603$         	&$3.841e-1$ 		&$0.470$\\
16	&$1.016e-1$        	& $0.733$        	&$2.399e-1$			&$ 0.537$\\
32	&$5.746e-2$        	& $0.823$        	&$1.606e-1$			&$0.579$\\
64	&$3.120e-2$        	& $0.881$       	&$1.056e-1$			&$0.606$\\
128	&$1.650e-2$        	& $0.919$      	&$6.850e-2$			&$0.624$\\
256	&$8.574e-3$        	& $0.944$         	&$4.406e-2$			&$0.637$\\
\hline
\end{tabular}
\vspace{0.2cm} 

Table 1. {\it Convergence study for mixed regime flows using mixed FEM in 2D.}
\end{center} 
We observe that by refining the time step and mesh size, the $L^2$-errors for the density  and  the $L^s$-error for momentum decrease
with a rate higher the theoretical rates of convergence from Theorem~\ref{Err-ful}.

 {\bf Example 2.}  We test the convergence of our method with the functions  $F_1(y)=y^{-1/2}+1+y$ and $F_2(y)=0.95 y^{-1/2}+1+0.95y$.
The analytical solution is computed according to $F_2(y)$ given as below 
\beqs
\rho(\x,t)=\frac 1 {\sqrt 2}\Big(0.95 e^{-t} +e^{-2t}+ 0.95 e^{-4t}\Big)(x_1+x_2), \, \text{ and }\,   
 \m(\x,t) =-e^{-2t}\left(\frac 1{\sqrt 2}, \frac 1{\sqrt 2}\right)  \quad \forall  (\x,t)\in \Omega\times I.
\eeqs
The forcing term $f$, initial condition and boundary condition accordingly are  
\beqs
f(\x,t)= -\frac 1{\sqrt 2}\Big(0.95 e^{-t} +2e^{-2t}+ 3.80 e^{-4t}\Big)(x_1+x_2),\quad \rho_0(x) =1.45\sqrt{2}(x_1+x_2), 
\eeqs
\beqs
 \psi(\x,t)=\frac 1 {\sqrt 2}\Big(0.95 e^{-t} +e^{-2t}+ 0.95 e^{-4t}\Big)\begin{cases}  
x_2  & \text { on } \{0\}\times[0,1] ,\\
1+x_2  & \text { on } \{1\}\times[0,1], \\
x_1  & \text { on } [0,1]\times \{0\}, \\
1+x_1  & \text { on } [0,1] \times \{1\}.
\end{cases}
\eeqs
We use $\|\rho_h^{(1)}-\rho_h^{(2)}\|_{L^2}$ and $\|\m_h^{(1)}-\m_h^{(2)}\|_{L^3}$ as the criterion to measure the dependence of solutions on the coefficients of $g$. The numerical results are listed  in Table 2.

\vspace{0.3cm}  
\begin{center}
\begin{tabular}{l| c| c| c| c}
\hline
&    &   &   \\
 $\mathcal N$    & \qquad $\norm{ \brho_h^{(1)}-\brho_h^{(2)} }_{L^2(\Omega)}$ \qquad  &  \quad Rates  \quad&  \quad$\norm{ \m_h^{(1)}-\m_h^{(2)} }_{L^3(\Omega)}$ \quad & \quad  Rates \quad  \\
&    &   &   \\
\hline
4		& $1.328e-2$      	& --         		& $1.126e-2$		&--\\
8		& $7.946e-3$       	& $0.741$         	& $7.225e-3$		&$0.640$\\
16		&  $4.423e-3$      	& $0.845$        	& $4.595e-3$		&$0.653$\\
32		&  $2.365e-3$     	& $0.903$       	& $2.912e-3$		&$0.658$\\
64		&  $1.236e-3$     	& $0.937$      	& $1.840e-3$		&$0.662$\\
128		&  $6.364e-4$    	& $0.957$     	& $1.160e-3$		&$0.665$\\
256  		&$3.247e-4$        	& $0.971$       	& $7.313e-4$		&$0.666$\\
\hline
\end{tabular}

\vspace{0.2cm}
Table 2. {\it Study the dependence of the solution of mixed regime flows using mixed FEM in 2D.}
\end{center}
We observe that by refining mesh size and the time step, the $L^2$-errors for the density  and  the $L^\alpha$-error for momentum decrease
with a rate higher the theoretical rates of convergence. This seems to indicate that the practical rate  of convergences are better than the theories ones from Theorem~\ref{Dep-ful}.
\def\cprime{$'$} \def\cprime{$'$} \def\cprime{$'$}
\def\cprime{$'$}

 \bibliographystyle{abbrv}

\begin{thebibliography}{10}

\bibitem{ABHI1}
E.~Aulisa, L.~Bloshanskaya, L.~Hoang, and A.~Ibragimov.
\newblock {Analysis of generalized {F}orchheimer flows of compressible fluids
  in porous media}.
\newblock {\em J. Math. Phys.}, 50(10):103102, 44, 2009.

\bibitem{BASAK77}
P.~BASAK.
\newblock Non-{D}arcy flow and its implications to seepage problems.
\newblock {\em J. Irrig. Drainage. Div.}, 103(4):459--473, 1977.

\bibitem{BearBook}
J.~Bear.
\newblock {\em Dynamics of Fluids in Porous Media}.
\newblock Dover, New York, 1972.

\bibitem{Bloshanskaya2017}
L.~Bloshanskaya, A.~Ibragimov, F.~Siddiqui, and M.~Y. Soliman.
\newblock Productivity index for {D}arcy and pre-/post-{D}arcy flow (analytical
  approach).
\newblock {\em Journal of Porous Media}, 20(9):769--786, 2017.

\bibitem{BPS02}
J.~H. Bramble, J.~E. Pasciak, and O.~Steinbach.
\newblock On the stability of the $l^2$-projection in $h^1(\omega)$.
\newblock {\em Mathematics of Computation}, 71(237):pp. 147--156, 2002.

\bibitem{BF91}
F.~Brezzi and M.~Fortin.
\newblock {\em Mixed and hybrid finite element methods}, volume~15 of {\em
  Springer Series in Computational Mathematics}.
\newblock Springer-Verlag, New York, 1991.

\bibitem{CHIK1}
E.~Celik, L.~Hoang, A.~Ibragimov, and T.~Kieu.
\newblock Fluid flows of mixed regimes in porous media.
\newblock {\em Journal of Mathematical Physics}, 58(2):023102, 2017.

\bibitem{CHK2}
E.~Celik, L.~Hoang, and T.~Kieu.
\newblock Doubly nonlinear parabolic equations for a general class of
  {F}orchheimer gas flows in porous media.
\newblock {\em Nonlinearity}, 31(8):3617--3650, 2018.

\bibitem{CHK1}
E.~Celik, L.~Hoang, and T.~Kieu.
\newblock Generalized forchheimer flows of isentropic gases.
\newblock {\em Journal of Mathematical Fluid Mechanics}, 20:83--115, 2018.

\bibitem{Ciarlet78}
P.~G. Ciarlet.
\newblock {\em The finite element method for elliptic problems}.
\newblock North-Holland Publishing Co., Amsterdam, 1978.
\newblock Studies in Mathematics and its Applications, Vol. 4.

\bibitem{Dudgeon85}
C.~Dudgeon.
\newblock {\em Non-{D}arcy flow of groundwater. {P}art 1. {T}heoretical,
  experimental and numerical studies. Report No. 162}.
\newblock Water research laboratory, {T}he {U}niversity of {N}ew {S}outh
  {W}ales, 1985.

\bibitem{ForchheimerBook}
P.~{F}orchheimer.
\newblock {\em Wasserbewegung durch Boden Zeit}, volume~45.
\newblock Ver. Deut. Ing., 1901.

\bibitem{M2AN75}
R.~Glowinski and A.~Marroco.
\newblock Sur l'approximation, par \'el\'ements finis d'ordre un, et la
  r\'esolution, par p\'enalisation-dualit\'e d'une classe de probl\`emes de
  dirichlet non lin\'eaires.
\newblock {\em ESAIM: Mathematical Modelling and Numerical Analysis -
  Mod\'elisation Math\'ematique et Analyse Num\'erique}, 9(R2):41--76, 1975.

\bibitem{HI1}
L.~Hoang and A.~Ibragimov.
\newblock Structural stability of generalized {F}orchheimer equations for
  compressible fluids in porous media.
\newblock {\em Nonlinearity}, 24(1):1--41, 2011.

\bibitem{HI2}
L.~Hoang and A.~Ibragimov.
\newblock {Qualitative Study of Generalized {F}orchheimer Flows with the Flux
  Boundary Condition}.
\newblock {\em Adv. Diff. Eq.}, 17(5--6):511--556, 2012.

\bibitem{HIKS1}
L.~Hoang, A.~Ibragimov, T.~Kieu, and Z.~Sobol.
\newblock Stability of solutions to generalized {F}orchheimer equations of any
  degree.
\newblock {\em J. Math. Sci.}, 210(4):476--544, 2015.

\bibitem{HIK1}
L.~T. Hoang, A.~Ibragimov, and T.~T. Kieu.
\newblock One-dimensional two-phase generalized {F}orchheimer flows of
  incompressible fluids.
\newblock {\em J. Math. Anal. Appl.}, 401(2):921--938, 2013.

\bibitem{HIK2}
L.~T. Hoang, A.~Ibragimov, and T.~T. Kieu.
\newblock A family of steady two-phase generalized {F}orchheimer flows and
  their linear stability analysis.
\newblock {\em J. Math. Phys.}, 55:123101, 2014.

\bibitem{HK1}
L.~T. Hoang and T.~T. Kieu.
\newblock Interior estimates for generalized forchheimer flows of slightly
  compressible fluids.
\newblock {\em Advanced Nonlinear Studies}, 17(4):739--767, 2017.

\bibitem{HK2}
L.~T. Hoang and T.~T. Kieu.
\newblock Global estimates for generalized {F}orchheimer flows of slightly
  compressible fluids.
\newblock {\em Journal d'Analyse Math\'ematique}, 137(1565-8538):1--55, 2019.

\bibitem{HKP1}
L.~T. Hoang, T.~T. Kieu, and T.~V. Phan.
\newblock {Properties of generalized {F}orchheimer flows in porous media}.
\newblock {\em J. Math. Sci.}, 202(2):259--332, 2014.

\bibitem{JS1995}
D.~Jerison and C.~Kenig.
\newblock The inhomogeneous dirichlet problem in lipschitz domains.
\newblock {\em Journal of Functional Analysis}, 130(1):161 -- 219, 1995.

\bibitem{JT81}
C.~Johnson and V.~Thom{\'e}e.
\newblock Error estimates for some mixed finite element methods for parabolic
  type problems.
\newblock {\em RAIRO Anal. Num\'er.}, 15(1):41--78, 1981.

\bibitem{K17}
T.~Kieu.
\newblock A mixed finite element approximation for {D}arcy–{F}orchheimer
  flows of slightly compressible fluids.
\newblock {\em Applied Numerical Mathematics}, 120:141 -- 164, 2017.

\bibitem{K2}
T.~Kieu.
\newblock Existence of a solution for generalized {F}orchheimer flow in porous
  media with minimal regularity conditions.
\newblock {\em Journal of Mathematical Physics}, 61(1):013507, 2020.

\bibitem{K3}
T.~Kieu.
\newblock Solution of the mixed formulation for generalized {F}orchheimer flows
  of isentropic gases.
\newblock {\em Journal of Mathematical Physics}, 61(8):081501, 2020.

\bibitem{K1}
T.~T. Kieu.
\newblock Analysis of expanded mixed finite element methods for the generalized
  forchheimer flows of slightly compressible fluids.
\newblock {\em Numerical Methods for Partial Differential Equations},
  32(1):60--85, 2016.

\bibitem{KIM199675}
M.-Y. Kim, F.~Milner, and E.-J. Park.
\newblock Some observations on mixed methods for fully nonlinear parabolic
  problems in divergence form.
\newblock {\em Applied Mathematics Letters}, 9(1):75 -- 81, 1996.

\bibitem{PG16}
P.~Knabner and G.~Summ.
\newblock {Solvability of the mixed formulation for {D}arcy-{F}orchheimer flow
  in porous media}.
\newblock {2017}.
\newblock manuscript.

\bibitem{KUNDU2016278}
P.~Kundu, V.~Kumar, and I.~M. Mishra.
\newblock Experimental and numerical investigation of fluid flow hydrodynamics
  in porous media: Characterization of pre-{D}arcy, {D}arcy and non-{D}arcy
  flow regimes.
\newblock {\em Powder Technology}, 303:278 -- 291, 2016.

\bibitem{MR0259693}
J.-L. Lions.
\newblock {\em Quelques m{\'e}thodes de r{\'e}solution des probl{\`e}mes aux
  limites non lin{\'e}aires}.
\newblock Dunod, 1969.

\bibitem{Muskatbook}
M.~Muskat.
\newblock {\em The flow of homogeneous fluids through porous media}.
\newblock McGraw-Hill Book Company, inc., 1937.

\bibitem{NieldBook}
D.~A. Nield and A.~Bejan.
\newblock {\em {Convection in porous media}}.
\newblock Springer-Verlag, New York, fourth edition, 2013.

\bibitem{EJP05}
E.-J. Park.
\newblock Mixed finite element methods for generalized {F}orchheimer flow in
  porous media.
\newblock {\em Numer. Methods Partial Differential Equations}, 21(2):213--228,
  2005.

\bibitem{RT06}
P.~Raviart and J.~Thomas.
\newblock {\em A Mixed Finite Element Method for Second Order Elliptic
  Problems}, volume 606, pages 292--315.
\newblock 11 2006.

\bibitem{s97}
R.~E. Showalter.
\newblock {\em Monotone operators in {B}anach space and nonlinear partial
  differential equations}, volume~49 of {\em Mathematical Surveys and
  Monographs}.
\newblock American Mathematical Society, Providence, RI, 1997.

\bibitem{SSHI2016}
F.~Siddiqui, M.~Y. Soliman, W.~House, and A.~Ibragimov.
\newblock Pre-{D}arcy flow revisited under experimental investigation.
\newblock {\em Journal of {A}nalytical {S}cience and {T}echnology}, 7(2), 2016.

\bibitem{SJ78}
J.~Simon.
\newblock Caractérisation d’espaces fonctionnels.
\newblock {\em Bollettino della Unione Matematica Italiana. Series V. B}, 01
  1978.

\bibitem{SONI1978231}
J.~Soni, N.~Islam, and P.~Basak.
\newblock An experimental evaluation of non-darcian flow in porous media.
\newblock {\em Journal of Hydrology}, 38(3):231 -- 241, 1978.

\bibitem{SoniIslamBasak78}
J.~Soni, N.~Islam, and P.~Basak.
\newblock An experimental evaluation of non-{D}arcian flow in porous media.
\newblock {\em Journal of {H}ydrology}, 38(3-4):231--241, 1978.

\bibitem{StraughanBook}
B.~Straughan.
\newblock {\em {Stability and wave motion in porous media}}, volume 165 of {\em
  {Applied Mathematical Sciences}}.
\newblock Springer, New York, 2008.

\bibitem{VazquezPorousBook}
J.~L. V{\'a}zquez.
\newblock {\em {The porous medium equation}}.
\newblock {Oxford Mathematical Monographs}. The Clarendon Press Oxford
  University Press, Oxford, 2007.
\newblock Mathematical theory.

\bibitem{Ward64}
J.~C. Ward.
\newblock {Turbulent flow in porous media.}
\newblock {\em Journal of the Hydraulics Division, Proc. Am. Soc. Civ. Eng.},
  90(HY5):1--12, 1964.

\bibitem{z90}
E.~Zeidler.
\newblock {\em Nonlinear functional analysis and its applications. {II}/{B}}.
\newblock Springer-Verlag, New York, 1990.
\newblock Nonlinear monotone operators, Translated from the German by the
  author and Leo F. Boron.

\end{thebibliography}
 
 \end{document}